%% file: TACA_arXiv.tex
\documentclass[12pt, english]{article}
\usepackage{babel}
\usepackage[latin1]{inputenc}
\usepackage{amsthm}
\usepackage{amsfonts,amsmath,amstext,amsbsy,euscript,amssymb, graphics}
\newcommand{\onematch}[1]{
  \begin{scope}[rotate=#1]
     \link;
  \end{scope}
}

\newcommand{\lastmatch}[1]{
  \begin{scope}[rotate=#1]
     \match;
  \end{scope}
}

\theoremstyle{plain}
\newtheorem{thm}{Theorem}
\newtheorem{cor}{Corollary}

\newtheorem{lemma}[thm]{Lemma}

\usepackage{pgf}
\usepackage{tikz}
\usetikzlibrary{arrows}
\usetikzlibrary{decorations.pathmorphing}
\usetikzlibrary{backgrounds}
\usetikzlibrary{positioning}
\usetikzlibrary{fit}
\usepackage{graphicx}
\usepackage{standalone}
\title{Impartial games emulating one-dimensional cellular automata and undecidability}
\author{Urban Larsson\\ Department of Mathematical Sciences\\
Chalmers University of Technology and G\"oteborg University}
\date{\today}

\begin{document}
\maketitle

\begin{abstract}
We study two-player \emph{take-away} games whose outcomes emulate two-state one-dimensional cellular automata, such as Wolfram's rules 60 and 110. Given an initial string consisting of a central data pattern and periodic left and right patterns, the rule 110 cellular automaton was recently proved Turing-complete by Matthew Cook. Hence, many questions regarding its behavior are algorithmically undecidable. We show that similar questions are undecidable for our \emph{rule 110} game.
\end{abstract}

\colorlet{brownblack}{black!80!brown}
\colorlet{brownbl}{black!20!brown}
\colorlet{yellowhite}{yellow!50!white}
\colorlet{lightgray}{gray!70!white}
\def\link{ \draw[gray, fill = yellow] (0,0) rectangle (4.5, .2); \filldraw [fill=brownblack] (0,0.1) circle (4pt);}
\def\match{ \draw[lightgray, fill = yellowhite] (0,0) rectangle (4.5, .2); \filldraw [fill=brownbl] (0,0.1) circle (4pt);}

\section{Introduction}\label{S:1}
We study the inter-connections between two popular areas of mathematics, two-player combinatorial games e.g. \cite{BCG04} and cellular automata (CAs) \cite{N66,HU79,W84a,W84b,W84c, W86, W02}. We present an infinite class of games and prove that their \emph{outcomes} (or winning strategies) emulate corresponding one-dimensional CAs. In particular we study some recent results of Matthew Cook, \cite{C04,C08}, concerning algorithmic undecidability of Stephen Wolfram's well known \emph{elementary cellular automaton}, rule 110, and interpret these results in the setting of our games. The universality of the rule 110 automaton was conjectured by S. Wolfram in 1985. It is also discussed in the remarkable book, \cite{W02}.

Our games are played between two players and are purely combinatorial---there is no element of chance and no hidden information. They are similar to the \emph{take away} games found in \cite{G66,S70,Z96}. In such games the players take turns in removing tokens (coins, matches, stones) from a finite number of heaps, each with a given finite number of tokens. A \emph{ruleset} gives the legal moves of a game. The \emph{ending condition} is given by: a player who is not able to move loses and the other player wins. This is \emph{normal play}. In \emph{Mis\`ere play} the ending condition is reversed, a player who is not able to move wins. If the set of options does not depend on which player is about to move, then the game is \emph{impartial}, otherwise the game is \emph{partizan}. In this paper we study normal play impartial games. The outcome classes of these games are denoted P (previous player win) and N (next player win). That is, a position (game) is P if and only if the player whose turn it is to move loses. 
This gives a recursive characterization of the outcomes of all starting positions of a game and (unless there are drawn positions such as infinite loops where no player can force a win) there will be a partitioning of the set of game positions into N and P. 

For certain games the pattern of the two sets of outcomes are reasonably easy to understand. For example it is known that the set of P-positions of a one heap subtraction game, e.g \cite{BCG04}, with a finite number of moves---such as a heap of a finite number of tokens and the ruleset, remove one, two or five tokens---is eventually periodic. On the other hand, in \cite{L, LW} simple rulesets are studied which give rise to very complex pattern of P-positions. In the latter paper the classical one heap subtraction games are generalized to several heaps and, by emulating binary one-dimensional cellular automata with finite update functions, it is shown that for finite ``rulesets'' it is undecidable whether or not two games have the same sets of P-positions. 
It appears that links between CAs and two-player combinatorial games are uncommon in the literature, the only sources except \cite{LW} that we have found so far are \cite{F02,F12}. (See also \cite{DH, DH09} for interesting surveys on algorithms, complexity and combinatorial games.)

The cellular automata use simple rules for updating some discrete structure in discrete (time) steps. For the one-dimensional case we take a doubly infinite binary string as input. As mentioned in the previous paragraph, if we fix some simple initial string (such as a single ``1'' among ``0''s), it is well known that e.g. \cite{W02}, for finite update functions as the ``program'', many essential problems regarding the CAs behavior are algorithmically undecidable. In \cite{C04, C08}, a particularly simple instance of an update function is studied, Wolfram's rule 110: the value of any given cell remains ``0'' if the neighboring cell to the left is also a ``0''. It remains a ``1'' if at least one of the two nearest neighbors contains a ``0''. Otherwise the value switches. If the initial string is arbitrary (or in a sense too simple or too complicated \cite{N66}), for a fixed update function, questions of decidability do not appear interesting. It turns out that a natural setting for the rule 110 CA is to code the ``program'' in a central finite part of the initial string together with certain periodic left and right patterns. Under these assumptions it is shown, in \cite{C04,C08}, that many questions regarding the rule 110 CA are undecidable. In fact, universal Turing machines with the least known number of states and symbols were constructed by simulating the rule 110 automaton. In contrast the simpler Wolfram's rules 60, defined in Section \ref{S1.1}, and 90 are known to be decidable under the above assumptions.

Since the one-dimensional cellular automata generate two-dimensional patterns over ``time'', it suffices to let our games be played on 2 finite heaps, a \emph{time-heap} and a \emph{tape-heap}, a terminology introduced in \cite{LW} but where several heaps were used in the simulation of cellular automata. Our variation of take away games belong to a different class than the classical subtraction games, namely the number of tokens a player is allowed to remove depends in some way on the previous player's move (so that in fact our game positions will be represented by ordered triples of non-negative integers). Such games are sometimes called \emph{move-size dynamic} e.g. \cite{L09} and in fact, we will adapt an idea from that paper: the move options on one of the heaps will depend on the other player's move on the other heap, although the context is quite different here (in that paper the context is to produce a game with the same outcome as the classical game of Wythoff Nim).

As an introductory example we begin by studying a game which emulates Wolfram's rule 60 CA and in a particularly simple setting, namely where the outcomes form Pascal's triangle modulo 2. In Section \ref{S2} we define a large class of games and prove that they emulate a certain family of cellular automata. Then, in Section \ref{S:T} we introduce simpler impartial ``triangle placing'' games with equivalent outcomes as those in Section \ref{S2}. In Section \ref{S3} we discuss how some undecidability problems for the rule 110 cellular automaton are interpreted in our setting of combinatorial games. 

\subsection{The rule 60 game and Pascal's triangle}\label{S1.1}
A position of the rule 60 game consists of a finite heap of matches and a finite heap of tokens. The matches simulate ``time''. There are 2 players who alternate removal of tokens and matches according to the following rules. A move consists of two parts: (1) at least one match is removed, at most the whole heap; (2) at most $m_p$ tokens are removed (possibly none), where $m_p$ denotes the number of matches removed by the other player in the previous move. A player may not remove the whole heap of matches, as in (1), unless all tokens are also removed, as in (2). A player who is unable to move loses. The other player wins. 
See Figures \ref{F0} and \ref{F01}.

\begin{figure}[ht]
\centering
\input{matches3.tex}\caption{The previous player removed the right-most match. Hence at most one token may be removed, which means that no move is possible and hence the previous player wins.}\label{F0}
\vspace{1 cm}
\input{matches4.tex}\caption{In this game, the next player wins by removing the last match together with both tokens.}\label{F01}
\end{figure}

The update rule of Wolfram's rule 60 cellular automaton (CA) is as follows: Assign arbitrary binary values to $a_x^0$ for all integers $x$. For $y>0$, let $a_x^y=0$ if $a_{x-1}^{y-1}=a_x^{y-1}$ and otherwise let $a_x^y=1$. In other words $a_x^y=f(a_{x-1}^{y-1},a_x^{y-1})$, where $f(i,j)=i\oplus j$, the operation being binary addition without carry, the ``Xor'' gate. The two-dimensional patterns obtained by this cellular automaton are algorithmically decidable given that the initial one-dimensional pattern is sufficiently simple, say left and right periodic together with a central data ``program''. In particular, if a spatial pattern consists in a single $1$, say $a_1^i=1$ if and only if $i=1$, then the updates correspond precisely to Pascal's triangle modulo~2. In fact, the outcomes of the rule 60 game correspond to the updates of the rule 60 CA with an initial string of the form $\ldots 000111\ldots$.

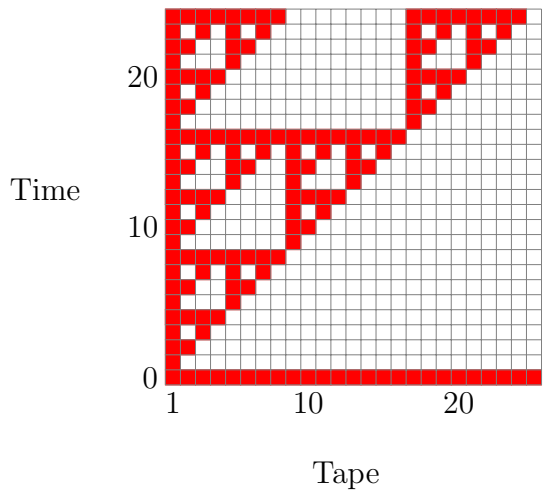
\begin{figure} [h]
\begin{center}
\begin{tikzpicture} [scale = 0.2]

\foreach \x/\y in {0/0, 1/0, 2/0, 3/0, 4/0, 5/0, 6/0, 7/0, 8/0, 9/0, 10/0, 11/0, 12/0, 13/0, 14/0, 15/0, 16/0, 17/0, 18/0, 19/0, 20/0, 21/0, 22/0, 23/0, 24/0, 0/1,0/2,0/3,0/4,0/5,0/6,0/7,0/8,0/9,0/10,0/11,0/12,0/13,0/14,0/15,0/16,0/17,0/18,0/19,0/20,0/21,0/22,0/23,0/24,
1/2,
2/3,
0/4, 1/4, 2/4, 3/4,
4/5,
1/6, 4/6, 5/6,
2/7, 4/7, 6/7,
1/8, 2/8, 3/8, 4/8, 5/8, 6/8, 7/8,
8/9,
1/10,
2/11,
0/12, 1/12, 2/12, 3/12,
4/13,
1/14, 4/14, 5/14,
2/15, 4/15, 6/15,
1/16, 2/16, 3/16, 4/16, 5/16, 6/16, 7/16,
8/10, 9/10,
8/11, 10/11,
8/12, 9/12, 10/12, 11/12,
8/13, 12/13,
8/14, 9/14, 12/14, 13/14,
8/15, 10/15, 12/15, 14/15,
8/16, 9/16, 10/16, 11/16, 12/16, 13/16, 14/16, 15/16,
16/17,
16/18, 17/18,
16/19, 18/19,
16/20, 17/20, 18/20, 19/20,
16/21, 20/21,
16/22, 17/22, 20/22, 21/22,
16/23, 18/23, 20/23, 22/23,
16/24, 17/24, 18/24, 19/24, 20/24, 21/24, 22/24, 23/24,
1/18,
2/19,
1/20, 2/20, 3/20,
4/21,
1/22, 4/22, 5/22,
2/23, 4/23, 6/23,
1/24, 2/24, 3/24, 4/24, 5/24, 6/24, 7/24
}
 \filldraw[color = red] (\x, \y) rectangle (\x+1, \y+1);

\draw[step=1cm,gray,very thin] (0, 0) grid (25,25); 

\draw (0.5 ,-1.2) node {$1$};
\draw (9.5, -1.2) node {$10$};
\draw (19.5, -1.2) node {$20$};

\draw (-1, 0.5) node {$0$};
\draw (-1.5, 10.5) node {$10$};
\draw (-1.5, 20.5) node {$20$};

\draw (12, -6) node {Tape};
\draw (-8, 13) node {Time};

\end{tikzpicture}
\end{center}

\caption{The CA given by $f(x,y) = x\oplus y$ (Wolfram's 
rule 60) together with an initial string of the form $\ldots 0011\dots $, the ``1''s correspond to red cells and the initial ``0''s are omitted. (Note that time flows upwards). By Theorem \ref{T:rule60}, for example the positions $(4,5,3), (8,9,7),(16,17,15) \ldots $ are all P, whereas $(x,x,m_p)$ is N, for all $x$ and $m_p$.}\label{F1}
\end{figure}

\begin{thm}\label{T:rule60}
Let the initial condition of the rule 60 CA be $a_i^0 = 1$ if and only if $i > 0$. Then, a position in the rule 60 game with $x$ tokens, $y$ matches and where the previous player removed $m_p$ matches is a second player win if and only if $a_{x}^{y} = \ldots = a_{x}^{y+m_p-1} = 0$ and if $y > 0$ then $a_{x}^{y-1} = 1$.
\end{thm}

Since $a_i^1=1$ if and only if $i=1$, the outcomes of this game are given by Pascal's triangle modulo~2. See Figure \ref{F1}. This result follows from the main theorem in the next section.

\section{Take away games and cellular automata}\label{S2}

\begin{figure}[ht!]
\begin{center}
\includegraphics[width=0.49\textwidth]{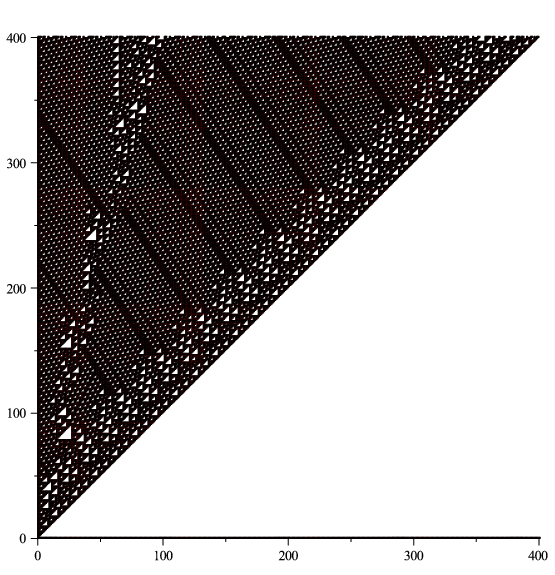}
\includegraphics[width=0.49\textwidth]{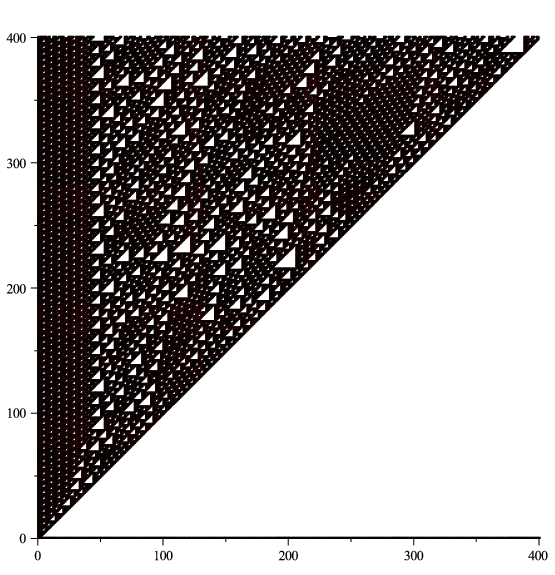}\vspace{0.3 cm}
\includegraphics[width=0.49\textwidth]{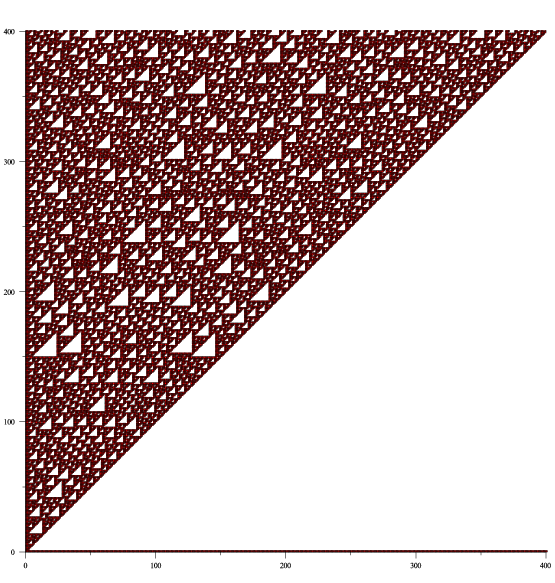}
\includegraphics[width=0.49\textwidth]{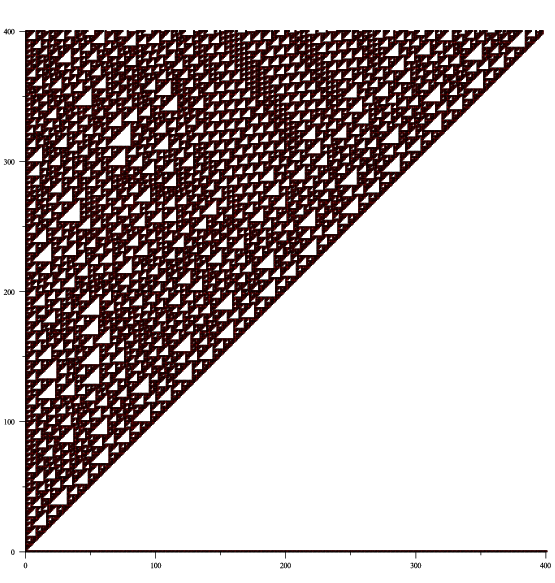}
\end{center}\caption{The updates of CA$(A, 1 , 0)$, CA$(A, 1 , 2)$, CA$(A, 1 , 1)$  and CA$(A, 3 , 0)$, for $A=\ldots0011\ldots$ and $y\le 400$. The first figure, which is rule 110, is a ``Class 4'' CA corresponding to Wolfram's classification of one-dimensional CAs. The third appears to be ``Class 3'', whereas the second and fourth are yet unclear. (As before, we let time flow upwards.)}\label{Feps} 
\end{figure}

In this section we define a generalization of the rule 60 game where the move options depend on two non-negative integer parameters, $\gamma $ and $\Gamma$, and the ending condition on a given ``black or white coloring'' of each token. Since our take away games are move-size dynamic the move options for the next player will depend in some precise manner, generalizing the rule 60 game, on the previous player's move. The idea is, roughly, that the games will emulate a special family of one-dimensional cellular automata for which, the update function will produce a zero in a given cell at a given positive time either if it reads two consecutive zeros (as in the rule 60 update) or if it reads a finite number of consecutive ones (generalizing the rule 60 update). We will return to the precise definition of the CAs, but a remark is in order here: by letting the games emulate the CAs as described, the pattern of ``discrete'' \emph{right isosceles triangles} as seen in Figure \ref{F1} will reappear in all of our games (see Figure \ref{Feps}). This theme will be further developed in Section \ref{S3}. 

A game position consists of a heap of $Y$ (unordered) matches and a heap of $X$ ordered and ``colored'' tokens, where $X$ and $Y$ are non-negative integers. Let $\tau_1\ldots \tau_X$ denote a finite binary string. Then the $i$th token is ``black'' if and only if $\tau_i = 1$, the $1^{st}$ token is at the bottom of the heap and the $X^{th}$ token is at the top. (It may be convenient to think of `non-positive tokens' as ``white'', but we will make our definitions independent of this.) The two players alternate turns in which they remove tokens and matches according to the following rules. A move consists of two parts:
\begin{itemize}
\item [(I)] remove first $t$ tokens from the top of the heap, where $0\le \gamma (m-1) \le t \le \gamma m + m_p + \Gamma $ and where $m_p$ denotes the number of matches removed by the other player in the previous move, with the exception that if the number of remaining tokens is less than $\gamma (m-1)$ then all of them can be removed, 
\item [(II)] remove $m$ matches, at least one and at most the whole heap. Given the removal of tokens as in (I), the whole heap of matches, $y$ of them, can be removed if and only if there is no black token among the top $y$ tokens. 
\end{itemize}

We connect the outcomes of these games to the patterns of cellular automata. Hence let us define the CAs of interest (see also Figure  \ref{Feps} for some examples). Assign an initial string of binary values $A = (a_x) = (a^0_x)$ for all integers $x$. The update rule of the cellular automaton CA$(A, \gamma , \Gamma)$ is as follows:  For $y>0$, let $a_x^y=0$ if $$a_{x-1}^{y-1}=a_i^{y-1}=0$$ or if $$a_{x-\Gamma-1}^{y-1}=\ldots =a_{x+\gamma}^{y-1}=1$$ and $a_x^y=1$ otherwise. Later we will let the variables $x$ and $y$ refer to the coordinates of a corresponding CA. 

Given game constants $\Gamma , \gamma$ and a doubly infinite binary string $A$, we denote a game family by G$(A,\gamma,\Gamma)$ and a specific game by G$(A,\gamma,\Gamma)^{(X,Y,m_p)}_{\xi}$, where $\xi$ together with the string $A$, determine the specific coloring of the game and $(X,Y,m_p)$ specifies the position. Precisely, the $X$ tokens are colored by the finite binary string $a_{\xi+1}\ldots a_{\xi +X}\subset A$ via the rule: $a_{\xi +i}=1$ iff the $i^{th}$ token is ``black'', that is $\tau_i=a_{\xi+i}$, for all $i$. The bottom token is colored according to the value of $a_{\xi +1}=\tau_1$ and the top token according to $a_{\xi +X}=\tau_X$.

If $A=\underline{0}$ then the first player wins (independent of the other variables), where $\underline{\cdot}$ denotes a periodic given pattern (infinite or doubly infinite). As we have seen in Section \ref{S:1}, if $A=\underline{0}\, \underline{1}$, precisely $a_x=1$ if and only if $x\ge 1$, $\xi=0$ and $\gamma = \Gamma = 0$, the CA describes the outcomes of the rule 60 game. In the next section we study the \emph{rule 110 game} which corresponds to $\gamma = 1$ and $\Gamma = 0$, see also Figure \ref{F:game3} for a particular position. Its outcome is illustrated in Figure \ref{F:8}, which leads us to the main result of this section.
\vspace{0.5 cm}
\begin{figure}[ht]
\centering
\input{matches1.tex}
\end{figure}
\vspace{0.5 cm}
\begin{thm}\label{T:2}
Let $A = (a_i)$ denote an initial condition of the cellular automaton CA$(A,\gamma , \Gamma)$ and let the game be G$(A, \gamma, \Gamma)_0^{(x+\gamma y, y, m_p)}$. If in addition $x\ge (\Gamma +1)y + m_p$, then the following conditions are equivalent:
\begin{itemize}
\item[(i)] The updates of the CA satisfy $a_{x}^{y} = \ldots = a_{x}^{y+m_p-1} = 0$ and if $y>0$ then $a_{x}^{y-1} = 1$. 
\item[(ii)] The position $(X,Y,m_p) = (x + \gamma y, y, m_p)$ is P, a previous player win. 
\end{itemize}
The same result holds in full generality with the initial condition of the CA exchanged for $A = \ldots 00a_1 a_{2}\ldots$. In particular this result holds whenever $-\gamma y\le x<(\Gamma +1)y+m_p $.
\end{thm}

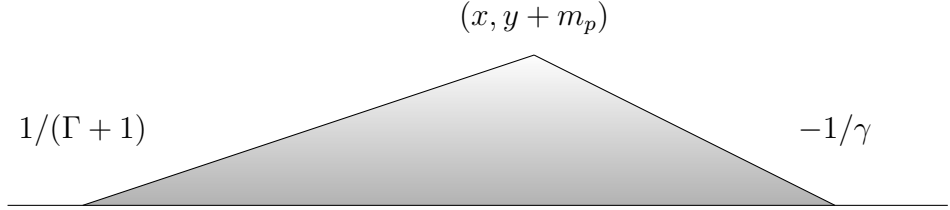
\begin{figure} [h]
\hspace{0.5 cm}
\begin{center}
\begin{tikzpicture} [scale = 0.5]
\draw (12, 2) node {$-1/\gamma$};
\draw (-8, 2) node {$1/(\Gamma +1)$};
\shade[top color=gray!1, bottom color=gray!60] (-8,0) -- (4,4)--(12,0)--cycle;
\draw (4, 5) node {$(x,y+m_p)$};
\draw (12,0) -- (4,4);
\draw (-8,0) -- (4,4);
\draw (-10,0) -- (15,0);
\end{tikzpicture}\caption{The value of the game position $(x+\gamma y,y,m_p)$ is influenced by the CA-cells bounded above by lines of slopes $1/(\Gamma +1)$ and  $-1/\gamma$ respectively.}\label{F:5}
\end{center}
\end{figure}

\begin{figure} [h]
\hspace{0.5 cm}
\begin{center}
\begin{tikzpicture} [scale = 0.5]
\draw (12, 2) node {$1/(\Gamma +1)$};
\draw (-8, 2) node {$-1/\gamma$};
\draw (0, -1) node {$a_i$};
\shade[top color=gray!1, bottom color=gray!90] (12,4)--(0,0)--(-8,4)--cycle;
\draw (0,0) -- (12,4);
\draw (0,0) -- (-8,4);
\draw (-10,0) -- (15,0);
\end{tikzpicture}\caption{The time-wise influence of the initial CA value $a_i$ is bounded below by lines of slopes  $-1/\gamma$ and $1/(\Gamma +1)$ respectively.}\label{F:6}
\end{center}
\end{figure}
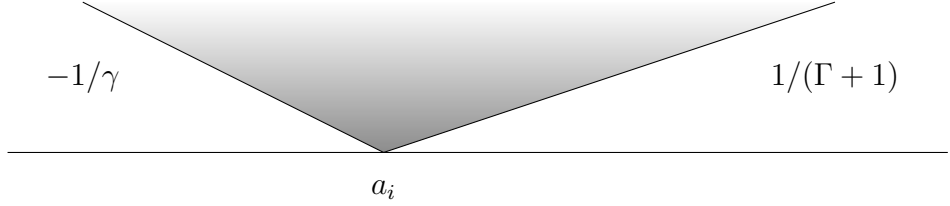

\noindent{\bf Proof.}
We begin by proving that if $x \ge (\Gamma + 1) y + m_p$, then the winning condition will depend on the coloring of the tokens and not of the empty tape-heap, that is after removal of tokens as in (I) and (the final) $y$ matches as in (II), there will remain at least $y$ tokens in the tape-heap. In this way we can guarantee that the outcomes are independent of the non-positive part of the CAs initial condition, see also Figures \ref{F:5} and \ref{F:6}. If the next player wishes to remove $m$ matches and keep the x-coordinate of the CA constant then a precise removal of $\gamma m$ tokens is required.  Therefore define $s=t-\gamma m$, where, as before, $t$ denotes the total number of tokens removed. Given $-\gamma \le s \le m_p + \Gamma$, we wish to show that $x - s \ge (\Gamma + 1)(y - m) + m = y+\Gamma (y-m)\ge y$, where the last inequality follows from $y\ge m$, since at most $y$ matches can be removed. We have that 
\begin{align*}
x-s&\ge  x-m_p-\Gamma\\ &\ge x+((\Gamma +1)y -x)-\Gamma\\&=(\Gamma +1)y-\Gamma\\&\ge (\Gamma +1)y-\Gamma m,
\end{align*}
where the last inequality follows from the fact that at least one match must be removed, which proves the claim.

For the rest of the proof we can ignore the finiteness of the heap of tokens and assume the appropriate initial condition of the CA. We need to show that, if a position as in (ii) satisfies (i), then none of its options does (denoted by ``P$\rightarrow$ N'') and that if a position as in (ii) does not satisfy (i) then one of its options does (``N$\rightarrow$ P''). Let us begin with the former case.\\

``P$\rightarrow $N'': Suppose first that $(x+\gamma y,y,m_p)$ is of the form in (i). Then we need to show that each of its options is not of this form. We may assume that $y>0$ since otherwise there is no option. An option is of the form 
\begin{align}\label{option}
(x+\gamma y- t,y-m, m) 
\end{align}
and where $0 \le \gamma (m-1) \le t \le \gamma m + m_p + \Gamma $. By the assumption we have that $a_{x}^{y-1}=1$. Also, item (i) together with the updates of the CA give $a_{x-i}^{y}=0$ for $0\le i\le m_p-1$. Altogether this gives that 
\begin{align}\label{Npos}
a_{x-j}^{y-1}=1 
\end{align}
for all $-\gamma \le j\le m_p+\Gamma$. By using (\ref{option}), we wish to derive $x'$ corresponding to the x-coordinate for the new CA cell. If it is bounded by $x-m_p-\Gamma \le x'\le x+\gamma$ then we are done with this part of the proof. Namely, the new y-coordinate of the same form as ``$y+m_p-1$'' equals $(y-m)+m-1=y-1$, which, by (\ref{Npos}), suffices to prove the claim. The y-coordinate of the next player's option is $y-m$, which gives $x'$ (defined above) as $x+\gamma y - t - \gamma (y-m) = x - t + \gamma m$. Thus the problem reduces to show that $-m_p-\Gamma \le -t+\gamma m\le \gamma$, that is that $t\le \gamma m + m_p+\Gamma$ and $\gamma (m-1)\le t$ which is true. \\

``$N\rightarrow P$'': For this case we have to show that it is possible to find an option of the form in (i) whenever one is playing from a position not of this form. Suppose first that $a_x^{y-1}=0$. Then there is a least $i\ge 1$ such that $a_x^{y-i}=0$ but $a_x^{y-i-1}\ne 0$. Remove $m=i$ matches and $\gamma m$ tokens (all tokens if $x < \gamma m$). Then the new position is of the correct form (since we assume that $a_{x-\gamma m}=0$ if  $x < \gamma m$). Otherwise we may assume that $a_x^{y+i}=1$ for some least $0\le i\le m_p-1$. By the updates of the CA (and minimality) this gives that the cell $a_{x-i}^{y} = 1$. We may assume that $a_x^{y-1} = 1$, therefore since $a_{x-i+1}^y=0$ the update rules of the CA force $a_{x-i-\Gamma -1}^{y-1}=0$. Hence we wish to remove tokens so that $x'=x-i-\Gamma -1$ is the x-coordinate of the CA-cell corresponding to the new position, now we can continue as in the first argument in this paragraph. 
Using the same argument as in the first part of the proof, we need to show that $-\gamma\le i+\Gamma +1 \le m_p + \Gamma $. The lower bound is clear since $i\ge 0$ and the upper bound, since $i+1\le m_p$.
\hfill $\Box$\\

In Theorem \ref{T:2}, for convenience, we have implicitly set $\xi = 0$, but one can easily deduce that it holds for all $\xi$ since the string $A$ can be translated arbitrarily preserving the same time-wise CA patterns (but at different spatial locations). The following corollary of Theorem \ref{T:2} supplies information about this and about the N-positions of our games. 

\begin{cor}\label{C:1}
A position $(x+\gamma y,y,m_p)$, of the game in Theorem \ref{T:2}, is N if and only if the corresponding CA updates satisfy one of the following: (a) $a_x^{y+i}=1$ for some $i\in \{0,\ldots , m_p-1\}$ or (b) $a_x^{y-1}=0$. 

Shift the indices in $A$ by $\xi$ steps so that in particular the content of the new 0-cell becomes that of the old $\xi$-cell, that is define $a'_{x-\xi} = a_{x}$ for all $x$. Then Theorem \ref{T:2} and the first paragraph of this corollary hold with $A$ exchanged for $(a'_i)$ and each $x$ exchanged for $x-\xi$. 
\end{cor}

\begin{figure}[ht]
\centering
\input{CAex3.tex}
\end{figure}
Two consequences of these results are the following ``periodicity lemma'' and the subsequent ``convergence lemma''. If the condition $x\ge (\Gamma +1)y + m_p$ in Theorem \ref{T:2} is satisfied, where the number of tokens in the tape-heap is $X=x+\gamma y\ge (\Gamma +\gamma + 1)y + m_p$ and the number of matches in the time-heap is $Y=y$, then we say that the tape-heap is \emph{super-critical}. See also Figure \ref{F:5}. The point of making this definition is that many results translate immediately between the two systems for super-critical tape-heaps. 

\begin{lemma}\label{L:3}
Let $A=(a_i)$ denote a doubly infinite binary string. Then the following items are equivalent. 
\begin{itemize}
\item The CA$(A,\gamma, \Gamma)$ has two-dimensional eventually periodic updates, that is, there are universal constants $\delta, \rho$ such that, for all but finitely many $x$ and $y$, $a_x^y = a_{x+\delta}^{y+\rho}$.
\item All games in G$(A,\gamma, \Gamma)$, with super-critical tape-heaps, have two-dimensional eventually periodic outcomes: that is, there are universal constants $\rho',\delta' $ such that, for all $\xi$ and for all $m_p$, for all but finitely many $X$ and $Y$, the outcomes of the positions $(X, Y, m_p)$ and $(X + \delta', Y + \rho', m_p)$, both with super-critical tape-heaps, are identical.
\end{itemize}
\end{lemma}

\noindent{\bf Proof.} For simplicity, by Corollary \ref{C:1}, we may assume $\xi=0$. Suppose that the CA has two-dimensional periodic patterns spatially and time wise with period $\delta$ and $\rho$ respectively. By Theorem \ref{T:2} we get that the outcomes of the positions $(x+\gamma y,y,m_p)$ and $(x+\delta + \gamma (y+\rho), y+\rho , m_p)$ are identical. Hence we can take $\delta'=\delta + \gamma \rho$.   

Suppose, on the other hand, that the outcomes of the positions $(X, Y, m_p) =  (x+\gamma y, y, m_p)$ and $(X + \delta', Y + \rho', m_p) = (x + \delta' + \gamma y, y + \rho' , m_p)$ are identical with super-critical tape-heaps. If they are both P, then, by Theorem \ref{T:2}, we have that $a_{x}^{y}=a_{x+\delta' -\rho'\gamma}^{y+\rho'}=0$. This implies that periodic N-positions of type (b) in Corollary \ref{C:1} can be dealt with analogously. Otherwise we get that $a_{x}^{y}=a_{x+\delta' -\rho'\gamma}^{y+\rho'} = 1$. In either case we can take $\delta = \delta'-\rho'\gamma$ and $\rho = \rho' $. Since we have assumed super-critical tape-heaps, the values of the CA correspond precisely to those of the games according to Theorem \ref{T:2}. 
\hfill $\Box$\\

The method in the proof actually says that the (three-dimensional) game positions define the pattern of the corresponding CA uniquely via its set of P-positions. This observation is used again in the next result. Let $A=(a_i)$ and $B=(b_i)$ denote doubly infinite binary strings. We say that the outcomes of the games in G$(A,\gamma,\Gamma)$ and G$(B,\gamma,\Gamma)$ \emph{converge} if, for all games on sufficiently large time-heaps with super-critical tape-heaps, for all $\xi$ and $m_p$, their outcomes are identical. The cellular automata CA$(A,\gamma,\Gamma)$ and CA$(B,\gamma,\Gamma)$ \emph{converge} if and only if, for all sufficiently large $y$, $a_x^y = b_x^y$ for all $x$. 

\begin{lemma}\label{L:4}
Let $A=(a_i)$ and $B=(b_i)$ denote doubly infinite binary strings. The outcomes of the games in G$(A,\gamma,\Gamma)$ and G$(B,\gamma,\Gamma)$ converge if and only if the cellular automata CA$(A,\gamma,\Gamma)$ and CA$(B,\gamma,\Gamma)$ converge.
\end{lemma}

\noindent{\bf Proof.} Since the tape-heaps are super-critical, by Theorem \ref{T:2}, the outcomes for the respective games in one of the families correspond precisely to the patterns of the corresponding CA. Hence, by Theorem \ref{T:2}, if the CAs converge, it follows that the outcomes of the game families converge. For the other direction we use a similar argument as in Lemma \ref{L:3}, the patterns of the CA is defined uniquely, given only the description of the P-positions of the game (via the move-size dynamic rule).
\hfill $\Box$

\section{A simpler impartial triangle placing game}\label{S:T}
Before we move on to the section of undecidability, let us discuss a simpler impartial game family which also emulates the class of cellular automata CA$(A,\gamma,\Gamma)$. It is simpler in the sense that it avoids the finite heaps condition of the take away game in the previous section. One may also argue that it is simpler in the sense that the move-size dependence is `built into' the position in a more efficient way, although the game rules are in essence two ways of saying the same thing. When we use the term ``triangle'' in this section we think of the set of discrete lattice points which are covered by a certain triangle shape with horizontal base and in case of an right isosceles triangle, the right angle is to the right.

\begin{figure} [h]
\hspace{0.5 cm}
\begin{center}
\begin{tikzpicture} [scale = 0.7]
\draw (5.6, 2) node {$-h/\gamma$};
\draw (-2.5, 2) node {$h/(h+\Gamma +1)$};
\fill[color=blue!80] (-2,0) -- (4,4)--(5,0)--cycle;
\fill[color=blue!30] (-1.6,0.3) -- (4,4)--(4.9,0.3)--cycle;
\fill[color=white!50] (0,0.3) -- (4,4)--(4,0.3)--cycle;
\draw (4, 4.8) node {$(x,y+h)$};
\draw (3.2,0.8) node {$(x,y)$};
\draw (5,0) -- (4,4);
\draw (-2,0) -- (4,4);
\draw (-2,0) -- (5,0);
\end{tikzpicture}\caption{The triangle position $(x,y,h)$, with its base-sensor $\{(x-\Gamma+1,y-1),\ldots(x+\gamma,y-1)\}$ in dark blue and its top at $(x,y+h)$. That is the slopes of the sides connecting the base to the top are  $h/(h+\Gamma +1)$ and $-h/\gamma$ respectively. Its IRT-sensor is a right isosceles triangle with the right angle at $(x,y)$ and height $h$.}\label{F:T}
\end{center}
\end{figure}
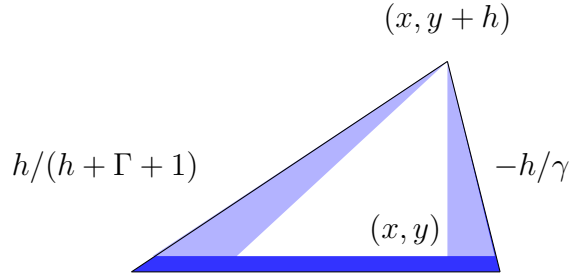

The upper half plane consists of all ordered pairs of integers $(x,y)$ with $y\ge 0$. The rules of the \emph{triangle placing game} T$(A,\gamma,\Gamma)$ are as follows. Let $(A,\gamma,\Gamma)$ be as in CA$(A,\gamma,\Gamma)$. Two players alternate in placing a \emph{$(\gamma, \Gamma)$-triangle} (triangle) of the shape in Figure \ref{F:T}, its size given by a non-negative integer $h$, on the upper half plane. Precisely, we denote a $(\gamma, \Gamma)$-triangle position by $(x,y,h)$, its \emph{top} covering the point $(x,y+h)$ and its \emph{base-sensor}, which has integer size $\Gamma +1 + h  + \gamma $, covering $\{(x-(\Gamma+1 +h),y-1),\ldots ,(x+\gamma, y-1)\}$. Here we only require that $y\ge 0$, so that it can be legal for a final triangle position to have the $(\gamma,\Gamma)$-triangle's base-sensor at the y-coordinate $-1$. The next player places another triangle, say $(x',y',h')$, of the same shape but possibly different size, with its top intersecting the base-sensor of the previous triangle, that is satisfying $x'\in \{x-(\Gamma+1 +h),\ldots ,x+\gamma \}$, with $y'+h'=y-1$ and $1\le h'\le y$.

As before, the ending condition is provided by the doubly infinite binary string $A = (a_x)$. Note that, by the rules of game the actual `game board' is bounded by a shape as in Figure \ref{F:5}. A player can place the $(\gamma,\Gamma)$-triangle $(x,0,h)$ if and only if $a_{x-h+1} + \ldots + a_x = 0$. Hence, as before, the ending condition does not depend on $\gamma$ and $\Gamma$. The \emph{IRT-sensor} of the triangle $(x,y,h)$ is another triangle, namely an \emph{isosceles right triangle} with its base covering the set $\{(x-h,y),\ldots ,(x,y)\}$ and the top covering the point $(x,y+h)$. By the update rule of the CA, another way of stating the ending condition is to require that if $y=0$ then the IRT-sensor covers only ``0''s in the underlying CA-cells. See also the second part of the proof of Theorem \ref{T:3}. 

\begin{thm}\label{T:3}
Given T$(A,\gamma,\Gamma)$, the triangle position $(x,y,h)$ is P if and only if its IRT-sensor covers only ``0''s and, if $y>0$, the \emph{base-sensor} covers only ``1''s in the update of CA$(A,\gamma,\Gamma)$, or equivalently $a_{x-h}^y + \ldots + a_x^y = 0$ and $a_{x-(\Gamma + 1 + h)}^{y-1} + \ldots + a_{x+\gamma }^{y-1} = \Gamma + h + \gamma + 2$.
\end{thm}

\noindent{\bf Proof.} For the ``$P\rightarrow N$'' direction, if $y=0$ we are done so suppose that $y>0$. If the IRT-sensor of the $(\gamma ,\Gamma)$-triangle $(x,y,h)$ covers only ``0''s and the base-sensor covers only ``1''s, then by the update of the CA, the next triangle position's IRT-sensor covers one of the ``1''s which was covered by the previous triangle position's base-sensor.

For the ``$N\rightarrow P$'' direction, suppose that the IRT-sensor of the triangle position $(x,y,h)$ covers a ``1''. Then the two CA-cells just below and to the left must cover at least one ``1''. Hence, viewing the ``1'', we can iterate this process until we approach the base of the IRT-sensor. Since this base covers a ``1'', by the ending condition this gives $y>0$, and, by the update rules of the CA, the base-sensor of the original triangle position $(x,y,h)$ must cover a ``0''. By the rules of triangle placing, the next player can use the ``0''-cell, at y-coordinate $y-1$, as the top of the next triangle and, by the update rule of the CA, choose the next $(\gamma,\Gamma)$-triangle suitably as to arrive at only ``0''s under its IRT-sensor and only ``1''s underneath its base-sensor if its y-coordinate is positive. 
\hfill $\Box$

\section{The rule 110 game and undecidability}\label{S3}
In this section we look into questions of decidability for the outcomes of our games. We focus at the take away games, but the results may equivalently be interpreted in the setting of the triangle placing games. By Lemma \ref{L:3} it is decidable whether the updates of any of our CAs eventually become two-dimensional periodic if and only if it is decidable whether the outcomes of the corresponding game do. By Lemma \ref{L:4} it is decidable whether two games converge if and only if it is decidable whether the corresponding CAs converge. As we discussed briefly in the introduction, questions of algorithmic decidability requires a finite input. For the CAs, this is achieved by letting the initial binary string be doubly periodic with a finite central data pattern. Such a binary string can encode the ending condition of a (family of) game(s), as described in previous sections. 

By recent results of Matthew Cook \cite{C04,C08} Wolfram's rule 110 CA, which in our notation is CA$(A,1,0)$, is particularly interesting, and hence also the game G$(A,1,0)$, which we also call the \emph{rule 110 game}. Let the initial binary string of this CA be of the form $A=\underline{L}C\underline{R}$, where $L, C$ and $R$ are finite binary strings (or equivalently integers coded in binary) and where, as before, $\underline{\cdot}$ denotes a periodic pattern. M. Cook proved the following results.

\begin{thm}[\cite{C04,C08}]
For finite binary strings $L$ and $R$ and a central finite data string $C$, it is algorithmically undecidable whether the rule 110 CA with $\underline{L}C\underline{R}$ as input ever produces a given binary string. 
\end{thm}

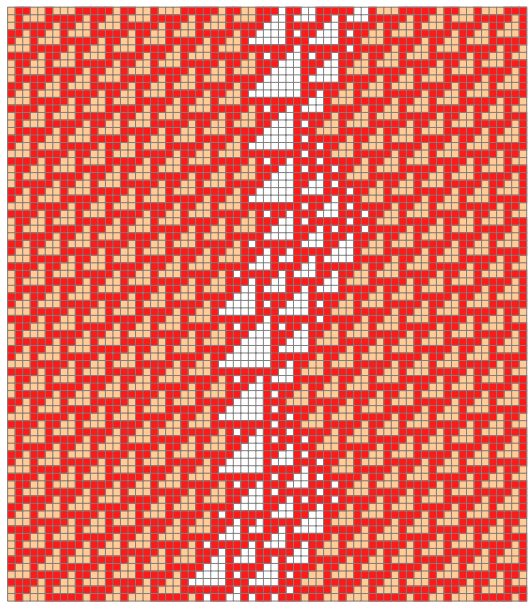
\begin{figure}[ht]
\centering
\input{CA1.tex}\caption{An F-glider embedded in Rule 110 ether.}\label{F:glider0}
\end{figure}

The proof uses that the particular binary string ``01101001101000'' is produced if and only if a certain ``F-glider'', see Figure \ref{F:glider0}, is created in the interaction of other ``gliders'' from the updates of the periodic $L$ pattern and the central $C$ pattern. Using \emph{cyclic tag-systems} \cite{C04,C08} a universal Turing machine is programmed to halt if and only if the given binary string occurs in the updates of the CA. This is how the rule 110 CA is proved undecidable. 
 One consequence of this result is that it is undecidable if the patterns in this CA will ultimately become two-dimensional periodic as defined in Lemma \ref{L:3}.

\begin{cor}[\cite{C04,C08}]
Let $L$, $R$ and $C$ denote finite binary strings. It is algorithmically undecidable whether (the central data pattern of) the rule 110 CA with $\underline{L}C\underline{R}$ as input are two-dimensional eventually periodic for a given $(space, time)$ period.
\end{cor}

An analogous corollary holds for our games. We do not need to assume super-critical tape-heaps although this is crucial in the construction since the result obviously depends on Theorem \ref{T:2} and Lemma \ref{L:3}.

\begin{cor}\label{C:2}
For fixed binary strings $L$, $R$ and a central data pattern $C$, it is algorithmically undecidable whether the outcomes of the rule 110 games, with tape-heaps, and thereby also ending conditions, defined by increasing finite sub-strings of $\underline{L}C\underline{R}$, are two-dimensional eventually periodic for a given $(space, time)$ period.
\end{cor}

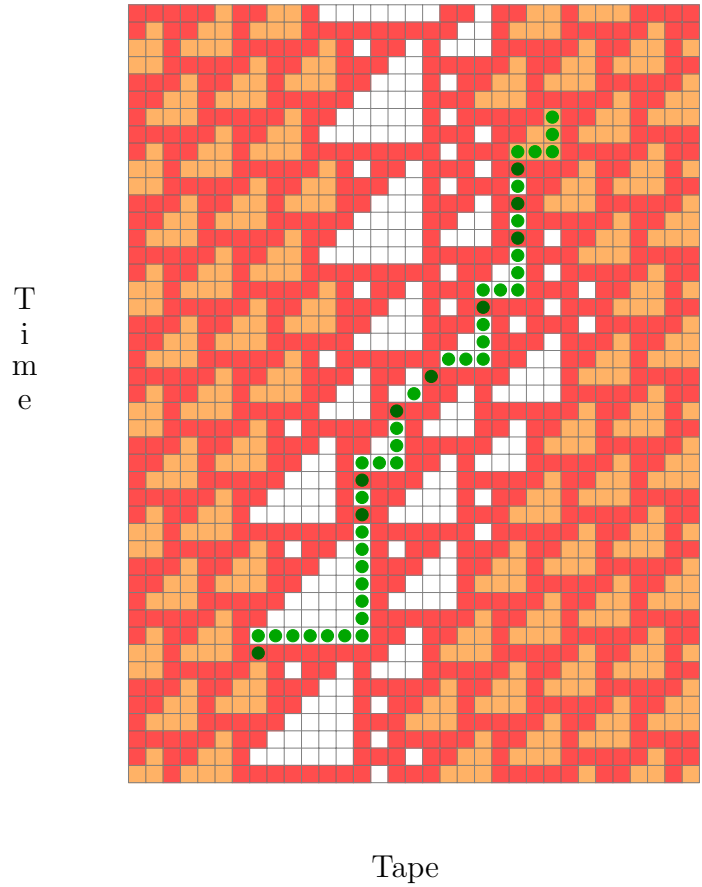
\begin{figure}[ht]
\centering
\input{CA2.tex}\caption{A good period of an F-glider embedded in Rule 110 ether. The green path of circles represents a sequence of alternating optimal moves traversing the F-glider. The light-green circles cover white CA-cells, building up reversed L shapes. Each such reversed L corresponds to a winning move, whereas the dark-green circles cover red cells and hence represent losing moves.}\label{F:glider}
\end{figure}

The halting problem for a universal Turing machine can be translated to the setting of our games via a finite ``optimal path'' (non-standard terminology) of alternating moves traversing the F-glider, illustrated in Figure \ref{F:glider}. An \emph{optimal path} of alternating moves is optimal for both players in the usual sense of \emph{perfect play}, that is no player can play any better, the \emph{winning player} can find a new P-position for each response of the \emph{losing player}. Using notation as in Corollary \ref{C:2}, we have the following result.

\begin{cor}
Let $L$, $R$ and $C$ denote finite binary strings. It is algorithmically undecidable whether a finite path of alternating moves is optimal in rule 110 games with $\underline{L}C\underline{R}$ ending conditions. 
\end{cor}

\noindent{\bf Proof.} The technical details of the proof are omitted. What is required is to show that the finite alternating path of moves in Figure \ref{F:glider} is optimal if and only if it traverses the F-glider. The light-green ``move-circles'' cover white ``CA-cells'' (0s) and the dark-green circles cover red cells (1s). The light-green ``reversed Ls'' correspond to winning moves: the number of vertical light-green (or equivalently horizontal) circles corresponds to $m_p$, the lower right corner to the heap sizes, the x-coordinate is the number of tokens minus the number of matches, the y-coordinate the number of matches, although the specific coordinates are not important here. Each dark-green circle corresponds to a move of the form, one match removed together with the top token. Hence, by comparing the game definitions with those of the CA and the F-glider \cite{C04}, the ``if'' direction is clear. 

For the `only if'-direction, by definition of this optimal path, the bottom dark-green move should be to an N position. This follows since the original CA-cell in the F-glider contains a ``1''. By Corollary \ref{C:1}, if it is of the form in (a) we are done. Hence suppose that it were of form (b), that is it covers a white CA-cell (and the cell just below is also white). This contradicts that the reversed light-green L just above is a winning move, namely either at least one of the CA cells in the base of the reversed L is red, or all CA-cells just below the base of the reversed L are white, in either case contradicting optimality of the move path. Let us next investigate the bottom reversed light-green L. Suppose that there was another underlying CA pattern for which this reversed L would still represent a P-position. Then the underlying white triangle would have needed to be at least as large (with the base at the same y-coordinate). But if it were any larger, by the update rule of the CA, the previous light-green move (represented by a single circle) would have covered a red cell, which again contradicts optimality. Similarly, if the dark-green circle on top of the bottom light-green reversed L would have been covering a white cell, it would still have been losing by Corollary \ref{C:1}, but so would the light green circle on top of it, hence destroying the alternating pattern of winning-losing moves. Thus, via analogous arguments the whole finite path is shown uniquely optimal for the F-glider. Since it is undecidable whether the F-glider is ever produced in the rule 110 CA,  with given $\underline{L}C\underline{R}$ input, it is also undecidable whether the alternating path of moves in Figure \ref{F:glider} ever becomes optimal in the corresponding rule 110 games. 
\hfill $\Box$ \\

Remark: our two imaginary players have become cooperative in proving Turing-completeness, they no longer compete in being the first to reach a given goal. In proving that a problem is undecidable, the competitive question regarding ``who wins?'' has been transferred to the question whether the players can cooperate in finding a certain optimal path of moves among games in the rule 110 family. Such a search could of course be carried out by our players for any other of our game families, but we do not yet know whether it would result in Turing-completeness. 

\section{Conclusion}\label{Conclusion}
The purpose of this paper has been to emulate well known cellular automata via impartial take-away games following the normal play convention. On the one hand, there is an `unintelligent system', in fact sometimes called zero player  game, with a very simple update function, which takes into account only the most recent history. On the other hand there are two combatants who battle intelligently in the attempt of being the first to reach a given goal; and where arbitrarily large moves were allowed. In spite the apparent big differences in the two systems we have showed that the `patterns' that the respective system produce correspond precisely and hence they are equivalent in many respect. 

Moreover, since one of the games emulates the rule 110 cellular automaton, we have demonstrated how the undecidability results from \cite{C04,C08} transfer to our setting. In the process we have discovered a simpler setting for an impartial triangle placing game. The corresponding undecidability results as presented in Section \ref{S3} for the games G$(A, \gamma,\Gamma)$ hold for the games T$(A,\gamma ,\Gamma)$. For example we can reformulate Corollary \ref{C:2} in this setting as follows.

\begin{cor}\label{C:T}
For fixed binary strings $L$, $R$ and a central data pattern $C$, it is algorithmically undecidable whether the P-positions of the game T$(\underline{L}C\underline{R}, 1, 0)$ ultimately become two-dimensional periodic.
\end{cor}

Some final remarks: returning to Lemma \ref{L:4}, and the first paragraph of Section \ref{S3}, it remains an open question whether convergence of rule 110 games is decidable given two $\underline{L}C\underline{R}$ ending conditions. Via private communication with Matthew Cook we understand that such results do not follow from the methods used in \cite{C04,C08}. Also, to our best knowledge, the problems of decidability discussed in this paper remain open for the CA$(\underline{L}C\underline{R},\gamma,\Gamma)$ and game families G$(\underline{L}C\underline{R},\gamma,\Gamma)$ for other combinations of $\gamma$ and $\Gamma$ than $\gamma \in\{0,1\}$ and $\Gamma = 0$. Many more CAs from \cite{W02} may have interesting interpretations as combinatorial games. A study of this territory may inspire new classifications of CAs. A further study of the class of CAs and games from this project will be given in \cite{L}.\\

\noindent{\bf Acknowledgments.} I thank David Wahlstedt and Matthew Cook for inspiring conversations on undecidability and Mike Weimerskirch for many useful comments.

\end{document}

%% file: matches3.tex

\begin{tikzpicture} [scale = 0.32]
\foreach \x/\y in {0/1,0/2}
 \filldraw[color = brownblack] (\x, \y) rectangle (\x+2.5, \y+1);
\draw[xstep=2.5cm,ystep=1cm, lightgray] (0, 1) grid (2.5,3); 
  \begin{scope}[xshift=6cm, yshift=2cm]
        \onematch{-10}
    \end{scope}
\begin{scope}[xshift=19cm,yshift=1.5cm]
  \begin{scope}[xshift=2.1cm,yshift=0.51cm]
        \lastmatch{150}
    \end{scope}
\end{scope}
\draw (-2,1) -- (12,1);
\end{tikzpicture}

%% file: matches4.tex

\begin{tikzpicture} [scale = 0.32]
\foreach \x/\y in {0/1,0/2}
 \filldraw[color = brownblack] (\x, \y) rectangle (\x+2.5, \y+1);
\draw[xstep=2.5cm,ystep=1cm, lightgray] (0, 1) grid (2.5,3); 
  \begin{scope}[xshift=6cm, yshift=2cm]
        \onematch{-10}
    \end{scope}
\begin{scope}[xshift=19cm,yshift=2cm]
  \begin{scope}[xshift=2.1cm,yshift=.41cm]
        \lastmatch{150}
    \end{scope}
 \begin{scope}[xshift=1.5cm,yshift=2.1cm]
        \lastmatch{220}
    \end{scope}
\end{scope}
\draw (-2,1) -- (12,1);
\end{tikzpicture}

%% file: matches1.tex



\colorlet{lightgray}{gray!70!white}

\begin{tikzpicture} [scale = 0.32]
\foreach \x/\y in {0/1,0/2,0/4}
 \filldraw[color = brownblack, opacity = .9] (\x, \y) rectangle (\x+2.5, \y+1);
\foreach \x/\y in {0/3,0/5,0/6}
\filldraw[color = white] (\x, \y) rectangle (\x+2.5, \y+1);
\draw[xstep=2.5cm, ystep=1cm, lightgray] (0, 1) grid (2.5,7); 
\end{tikzpicture}
\hspace{1.5cm}
\begin{tikzpicture} [scale = 0.32]
    \begin{scope}[xshift=0cm, yshift=1.31cm]
        \onematch{82}
    \end{scope}
  \begin{scope}[xshift=-1.2cm, yshift=4cm]
        \onematch{-10}
    \end{scope}

\begin{scope}[xshift=16cm, yshift=-0.5cm]
  \begin{scope}[xshift=2.1cm, yshift=2.31cm]
        \lastmatch{150}
    \end{scope}
  \begin{scope}[xshift=.2cm, yshift=2cm]
        \lastmatch{65}
    \end{scope}  
  \begin{scope}[xshift=1.4cm, yshift=4.1cm]
        \lastmatch{-170}
    \end{scope}
\end{scope}
\end{tikzpicture}
\caption{The rule 110 game position $(110100, 6, 2, 3)$ is a second player win.}\label{F:game3}

%% file: CAex3.tex
\begin{tikzpicture}[scale = 0.3]

\foreach \x/\y in {
0/0, 0/1, 
0/2, 0/3, 0/4, 0/5, 0/6, 0/7, 0/8, 
1/0, 1/1, 1/3, 1/4, 1/6, 1/8, 
2/1, 
2/4, 2/5, 2/6, 2/7, 2/8, 
3/0, 3/1, 3/5, 3/7, 3/8, 
4/1, 4/2, 4/3, 4/4, 4/5, 4/8, 
5/2, 5/4, 5/5, 
6/0, 6/1, 6/2, 6/5, 6/6, 6/7, 
7/0, 7/2, 7/6, 
8/0, 8/1, 8/2, 8/3, 8/4, 8/5, 
9/1, 9/3, 9/4, 
10/0, 10/1, 
11/0, 11/1, 
12/1
}
 \filldraw[color = red, opacity=0.9] (\x, \y) rectangle (\x+1, \y+1);

\foreach \x/\y in {
6/8, 
7/7,7/8, 
8/6, 8/7, 8/8,
9/5,9/6, 9/7,9/8,
10/4, 10/5, 10/6, 10/7, 10/8, 
11/3,11/4, 11/5, 11/6,11/7, 11/8, 
12/2, 12/3, 12/4, 12/5, 12/6,12/7,12/8, 
13/1,13/2, 13/3, 13/4, 13/5, 13/6, 13/7, 13/8, 
}
 \filldraw[color = gray, opacity=0.5] (\x, \y) rectangle (\x+1, \y+1);

\foreach \x/\y in {
0/0, 1/0,3/0,6/0,7/0,8/0,10/0,11/0
}
 \filldraw[color = black, opacity=0.6]   (\x, \y) rectangle (\x+1, \y+1);
\foreach \x/\y in {
3/2
}
 \filldraw[color = green, opacity=1]  (\x+.5, \y+.5) circle (11pt);
\foreach \x/\y in {
5/8
}
 \filldraw[color = blue, opacity=1] (\x+.5, \y+.5) circle (11pt);
\foreach \x/\y in {
0/1
}
 \filldraw[color = yellow, opacity=1]  (\x+.5, \y+.5) circle (11pt);
\foreach \x/\y in {
2/1
}
\colorlet{lightorange}{orange!50!yellow}
 \filldraw[color = lightorange, opacity=1]  (\x+.5, \y+.5) circle (11pt);
\draw[step=1cm, gray, ultra thin] (0, 0) grid (14,9); 

\draw (0.5 ,-1.2) node {$1$};
\draw (4.5, -1.2) node {$5$};
\draw (9.5, -1.2) node {$10$};

\draw (-1.5, 0.5) node {$0$};
\draw (-1.5, 5.5) node {$5$};


\end{tikzpicture}\caption{The filled circles indicate winning strategies of rule 110 game positions ($m_p$ omitted) for the finite ending condition $S=11010011101100$ together with CA updates. For example, the green position $(110100,6,2,m_p)$ is a second player win if and only if $1\le m_p\le 3$. For the other positions the first player wins independent of $m_p$. The cells in the gray area do not affect the winning strategy of the given positions.}\label{F:8}

%% file: CA1.tex
\begin{tikzpicture}[scale = 0.1]

\foreach \x/\y in {
10/22, 10/24, 10/25, 10/29, 10/31, 10/32, 10/36, 10/38, 10/39, 10/43, 10/45, 10/46, 10/50, 10/52, 10/53, 10/57, 10/59, 10/60, 10/64, 10/66, 10/67, 10/71, 10/73, 10/74, 10/78, 10/80, 10/81, 10/85, 10/87, 10/88, 10/92, 10/94, 10/95, 11/20, 11/21, 11/22, 11/25, 11/26, 11/27, 11/28, 11/29, 11/32, 11/33, 11/34, 11/35, 11/36, 11/39, 11/40, 11/41, 11/42, 11/43, 11/46, 11/47, 11/48, 11/49, 11/50, 11/53, 11/54, 11/55, 11/56, 11/57, 11/60, 11/61, 11/62, 11/63, 11/64, 11/67, 11/68, 11/69, 11/70, 11/71, 11/74, 11/75, 11/76, 11/77, 11/78, 11/81, 11/82, 11/83, 11/84, 11/85, 11/88, 11/89, 11/90, 11/91, 11/92, 11/95, 11/96, 11/97, 11/98, 12/21, 12/22, 12/26, 12/28, 12/29, 12/33, 12/35, 12/36, 12/40, 12/42, 12/43, 12/47, 12/49, 12/50, 12/54, 12/56, 12/57, 12/61, 12/63, 12/64, 12/68, 12/70, 12/71, 12/75, 12/77, 12/78, 12/82, 12/84, 12/85, 12/89, 12/91, 12/92, 12/96, 12/98, 13/22, 13/23, 13/24, 13/25, 13/26, 13/29, 13/30, 13/31, 13/32, 13/33, 13/36, 13/37, 13/38, 13/39, 13/40, 13/43, 13/44, 13/45, 13/46, 13/47, 13/50, 13/51, 13/52, 13/53, 13/54, 13/57, 13/58, 13/59, 13/60, 13/61, 13/64, 13/65, 13/66, 13/67, 13/68, 13/71, 13/72, 13/73, 13/74, 13/75, 13/78, 13/79, 13/80, 13/81, 13/82, 13/85, 13/86, 13/87, 13/88, 13/89, 13/92, 13/93, 13/94, 13/95, 13/96, 14/23, 14/25, 14/26, 14/30, 14/32, 14/33, 14/37, 14/39, 14/40, 14/44, 14/46, 14/47, 14/51, 14/53, 14/54, 14/58, 14/60, 14/61, 14/65, 14/67, 14/68, 14/72, 14/74, 14/75, 14/79, 14/81, 14/82, 14/86, 14/88, 14/89, 14/93, 14/95, 14/96, 15/20, 15/21, 15/22, 15/23, 15/26, 15/27, 15/28, 15/29, 15/30, 15/33, 15/34, 15/35, 15/36, 15/37, 15/40, 15/41, 15/42, 15/43, 15/44, 15/47, 15/48, 15/49, 15/50, 15/51, 15/54, 15/55, 15/56, 15/57, 15/58, 15/61, 15/62, 15/63, 15/64, 15/65, 15/68, 15/69, 15/70, 15/71, 15/72, 15/75, 15/76, 15/77, 15/78, 15/79, 15/82, 15/83, 15/84, 15/85, 15/86, 15/89, 15/90, 15/91, 15/92, 15/93, 15/96, 15/97, 15/98, 16/20, 16/22, 16/23, 16/27, 16/29, 16/30, 16/34, 16/36, 16/37, 16/41, 16/43, 16/44, 16/48, 16/50, 16/51, 16/55, 16/57, 16/58, 16/62, 16/64, 16/65, 16/69, 16/71, 16/72, 16/76, 16/78, 16/79, 16/83, 16/85, 16/86, 16/90, 16/92, 16/93, 16/97, 17/20, 17/23, 17/24, 17/25, 17/26, 17/27, 17/30, 17/31, 17/32, 17/33, 17/34, 17/37, 17/38, 17/39, 17/40, 17/41, 17/44, 17/45, 17/46, 17/47, 17/48, 17/51, 17/52, 17/53, 17/54, 17/55, 17/58, 17/59, 17/60, 17/61, 17/62, 17/65, 17/66, 17/67, 17/68, 17/69, 17/72, 17/73, 17/74, 17/75, 17/76, 17/79, 17/80, 17/81, 17/82, 17/83, 17/86, 17/87, 17/88, 17/89, 17/90, 17/93, 17/94, 17/95, 17/96, 17/97, 18/20, 18/24, 18/26, 18/27, 18/31, 18/33, 18/34, 18/38, 18/40, 18/41, 18/45, 18/47, 18/48, 18/52, 18/54, 18/55, 18/59, 18/61, 18/62, 18/66, 18/68, 18/69, 18/73, 18/75, 18/76, 18/80, 18/82, 18/83, 18/87, 18/89, 18/90, 18/94, 18/96, 18/97, 19/20, 19/21, 19/22, 19/23, 19/24, 19/27, 19/28, 19/29, 19/30, 19/31, 19/34, 19/35, 19/36, 19/37, 19/38, 19/41, 19/42, 19/43, 19/44, 19/45, 19/48, 19/49, 19/50, 19/51, 19/52, 19/55, 19/56, 19/57, 19/58, 19/59, 19/62, 19/63, 19/64, 19/65, 19/66, 19/69, 19/70, 19/71, 19/72, 19/73, 19/76, 19/77, 19/78, 19/79, 19/80, 19/83, 19/84, 19/85, 19/86, 19/87, 19/90, 19/91, 19/92, 19/93, 19/94, 19/97, 19/98, 20/21, 20/23, 20/24, 20/28, 20/30, 20/31, 20/35, 20/37, 20/38, 20/42, 20/44, 20/45, 20/49, 20/51, 20/52, 20/56, 20/58, 20/59, 20/63, 20/65, 20/66, 20/70, 20/72, 20/73, 20/77, 20/79, 20/80, 20/84, 20/86, 20/87, 20/91, 20/93, 20/94, 20/98, 21/20, 21/21, 21/24, 21/25, 21/26, 21/27, 21/28, 21/31, 21/32, 21/33, 21/34, 21/35, 21/38, 21/39, 21/40, 21/41, 21/42, 21/45, 21/46, 21/47, 21/48, 21/49, 21/52, 21/53, 21/54, 21/55, 21/56, 21/59, 21/60, 21/61, 21/62, 21/63, 21/66, 21/67, 21/68, 21/69, 21/70, 21/73, 21/74, 21/75, 21/76, 21/77, 21/80, 21/81, 21/82, 21/83, 21/84, 21/87, 21/88, 21/89, 21/90, 21/91, 21/94, 21/95, 21/96, 21/97, 21/98, 22/20, 22/21, 22/25, 22/27, 22/28, 22/32, 22/34, 22/35, 22/39, 22/41, 22/42, 22/46, 22/48, 22/49, 22/53, 22/55, 22/56, 22/60, 22/62, 22/63, 22/67, 22/69, 22/70, 22/74, 22/76, 22/77, 22/81, 22/83, 22/84, 22/88, 22/90, 22/91, 22/95, 22/97, 22/98, 23/21, 23/22, 23/23, 23/24, 23/25, 23/28, 23/29, 23/30, 23/31, 23/32, 23/35, 23/36, 23/37, 23/38, 23/39, 23/42, 23/43, 23/44, 23/45, 23/46, 23/49, 23/50, 23/51, 23/52, 23/53, 23/56, 23/57, 23/58, 23/59, 23/60, 23/63, 23/64, 23/65, 23/66, 23/67, 23/70, 23/71, 23/72, 23/73, 23/74, 23/77, 23/78, 23/79, 23/80, 23/81, 23/84, 23/85, 23/86, 23/87, 23/88, 23/91, 23/92, 23/93, 23/94, 23/95, 23/98, 24/22, 24/24, 24/25, 24/29, 24/31, 24/32, 24/36, 24/38, 24/39, 24/43, 24/45, 24/46, 24/50, 24/52, 24/53, 24/57, 24/59, 24/60, 24/64, 24/66, 24/67, 24/71, 24/73, 24/74, 24/78, 24/80, 24/81, 24/85, 24/87, 24/88, 24/92, 24/94, 24/95, 25/20, 25/21, 25/22, 25/25, 25/26, 25/27, 25/28, 25/29, 25/32, 25/33, 25/34, 25/35, 25/36, 25/39, 25/40, 25/41, 25/42, 25/43, 25/46, 25/47, 25/48, 25/49, 25/50, 25/53, 25/54, 25/55, 25/56, 25/57, 25/60, 25/61, 25/62, 25/63, 25/64, 25/67, 25/68, 25/69, 25/70, 25/71, 25/74, 25/75, 25/76, 25/77, 25/78, 25/81, 25/82, 25/83, 25/84, 25/85, 25/88, 25/89, 25/90, 25/91, 25/92, 25/95, 25/96, 25/97, 25/98, 26/21, 26/22, 26/26, 26/28, 26/29, 26/33, 26/35, 26/36, 26/40, 26/42, 26/43, 26/47, 26/49, 26/50, 26/54, 26/56, 26/57, 26/61, 26/63, 26/64, 26/68, 26/70, 26/71, 26/75, 26/77, 26/78, 26/82, 26/84, 26/85, 26/89, 26/91, 26/92, 26/96, 26/98, 27/22, 27/23, 27/24, 27/25, 27/26, 27/29, 27/30, 27/31, 27/32, 27/33, 27/36, 27/37, 27/38, 27/39, 27/40, 27/43, 27/44, 27/45, 27/46, 27/47, 27/50, 27/51, 27/52, 27/53, 27/54, 27/57, 27/58, 27/59, 27/60, 27/61, 27/64, 27/65, 27/66, 27/67, 27/68, 27/71, 27/72, 27/73, 27/74, 27/75, 27/78, 27/79, 27/80, 27/81, 27/82, 27/85, 27/86, 27/87, 27/88, 27/89, 27/92, 27/93, 27/94, 27/95, 27/96, 28/23, 28/25, 28/26, 28/30, 28/32, 28/33, 28/37, 28/39, 28/40, 28/44, 28/46, 28/47, 28/51, 28/53, 28/54, 28/58, 28/60, 28/61, 28/65, 28/67, 28/68, 28/72, 28/74, 28/75, 28/79, 28/81, 28/82, 28/86, 28/88, 28/89, 28/93, 28/95, 28/96, 29/20, 29/21, 29/22, 29/23, 29/26, 29/27, 29/28, 29/29, 29/30, 29/33, 29/34, 29/35, 29/36, 29/37, 29/40, 29/41, 29/42, 29/43, 29/44, 29/47, 29/48, 29/49, 29/50, 29/51, 29/54, 29/55, 29/56, 29/57, 29/58, 29/61, 29/62, 29/63, 29/64, 29/65, 29/68, 29/69, 29/70, 29/71, 29/72, 29/75, 29/76, 29/77, 29/78, 29/79, 29/82, 29/83, 29/84, 29/85, 29/86, 29/89, 29/90, 29/91, 29/92, 29/93, 29/96, 29/97, 29/98, 30/20, 30/22, 30/23, 30/27, 30/29, 30/30, 30/34, 30/36, 30/37, 30/41, 30/43, 30/44, 30/48, 30/50, 30/51, 30/55, 30/57, 30/58, 30/62, 30/64, 30/65, 30/69, 30/71, 30/72, 30/76, 30/78, 30/79, 30/83, 30/85, 30/86, 30/90, 30/92, 30/93, 30/97, 31/20, 31/23, 31/24, 31/25, 31/26, 31/27, 31/30, 31/31, 31/32, 31/33, 31/34, 31/37, 31/38, 31/39, 31/40, 31/41, 31/44, 31/45, 31/46, 31/47, 31/48, 31/51, 31/52, 31/53, 31/54, 31/55, 31/58, 31/59, 31/60, 31/61, 31/62, 31/65, 31/66, 31/67, 31/68, 31/69, 31/72, 31/73, 31/74, 31/75, 31/76, 31/79, 31/80, 31/81, 31/82, 31/83, 31/86, 31/87, 31/88, 31/89, 31/90, 31/93, 31/94, 31/95, 31/96, 31/97, 32/20, 32/24, 32/26, 32/27, 32/31, 32/33, 32/34, 32/38, 32/40, 32/41, 32/45, 32/47, 32/48, 32/52, 32/54, 32/55, 32/59, 32/61, 32/62, 32/66, 32/68, 32/69, 32/73, 32/75, 32/76, 32/80, 32/82, 32/83, 32/87, 32/89, 32/90, 32/94, 32/96, 32/97, 33/20, 33/21, 33/22, 33/23, 33/24, 33/27, 33/28, 33/29, 33/30, 33/31, 33/34, 33/35, 33/36, 33/37, 33/38, 33/41, 33/42, 33/43, 33/44, 33/45, 33/48, 33/49, 33/50, 33/51, 33/52, 33/55, 33/56, 33/57, 33/58, 33/59, 33/62, 33/63, 33/64, 33/65, 33/66, 33/69, 33/70, 33/71, 33/72, 33/73, 33/76, 33/77, 33/78, 33/79, 33/80, 33/83, 33/84, 33/85, 33/86, 33/87, 33/90, 33/91, 33/92, 33/93, 33/94, 33/97, 33/98, 34/21, 34/23, 34/24, 34/28, 34/30, 34/31, 34/35, 34/37, 34/38, 34/42, 34/44, 34/45, 34/49, 34/51, 34/52, 34/56, 34/58, 34/59, 34/63, 34/65, 34/66, 34/70, 34/72, 34/73, 34/77, 34/79, 34/80, 34/84, 34/86, 34/87, 34/91, 34/93, 34/94, 34/98, 35/20, 35/21, 35/24, 35/25, 35/26, 35/27, 35/28, 35/31, 35/32, 35/33, 35/34, 35/35, 35/38, 35/39, 35/40, 35/41, 35/42, 35/45, 35/46, 35/47, 35/48, 35/49, 35/52, 35/53, 35/54, 35/55, 35/56, 35/59, 35/60, 35/61, 35/62, 35/63, 35/66, 35/67, 35/68, 35/69, 35/70, 35/73, 35/74, 35/75, 35/76, 35/77, 35/80, 35/81, 35/82, 35/83, 35/84, 35/87, 35/88, 35/89, 35/90, 35/91, 35/94, 35/95, 35/96, 35/97, 35/98, 36/21, 36/25, 36/26, 36/28, 36/32, 36/34, 36/35, 36/39, 36/41, 36/42, 36/46, 36/48, 36/49, 36/53, 36/55, 36/56, 36/60, 36/62, 36/63, 36/67, 36/69, 36/70, 36/74, 36/76, 36/77, 36/81, 36/83, 36/84, 36/88, 36/90, 36/91, 36/95, 36/97, 36/98, 37/20, 37/21, 37/26, 37/27, 37/28, 37/29, 37/30, 37/31, 37/32, 37/35, 37/36, 37/37, 37/38, 37/39, 37/42, 37/43, 37/44, 37/45, 37/46, 37/49, 37/50, 37/51, 37/52, 37/53, 37/56, 37/57, 37/58, 37/59, 37/60, 37/63, 37/64, 37/65, 37/66, 37/67, 37/70, 37/71, 37/72, 37/73, 37/74, 37/77, 37/78, 37/79, 37/80, 37/81, 37/84, 37/85, 37/86, 37/87, 37/88, 37/91, 37/92, 37/93, 37/94, 37/95, 37/98, 38/20, 38/21, 38/27, 38/29, 38/30, 38/32, 38/36, 38/38, 38/39, 38/43, 38/45, 38/46, 38/50, 38/52, 38/53, 38/57, 38/59, 38/60, 38/64, 38/66, 38/67, 38/71, 38/73, 38/74, 38/78, 38/80, 38/81, 38/85, 38/87, 38/88, 38/92, 38/94, 38/95, 39/21, 39/22, 39/23, 39/24, 39/25, 39/26, 39/27, 39/30, 39/31, 39/32, 39/33, 39/34, 39/35, 39/36, 39/39, 39/40, 39/41, 39/42, 39/43, 39/46, 39/47, 39/48, 39/49, 39/50, 39/53, 39/54, 39/55, 39/56, 39/57, 39/60, 39/61, 39/62, 39/63, 39/64, 39/67, 39/68, 39/69, 39/70, 39/71, 39/74, 39/75, 39/76, 39/77, 39/78, 39/81, 39/82, 39/83, 39/84, 39/85, 39/88, 39/89, 39/90, 39/91, 39/92, 39/95, 39/96, 39/97, 39/98, 40/22, 40/24, 40/26, 40/27, 40/31, 40/33, 40/34, 40/36, 40/40, 40/41, 40/43, 40/47, 40/48, 40/50, 40/54, 40/55, 40/57, 40/61, 40/62, 40/64, 40/68, 40/70, 40/71, 40/75, 40/77, 40/78, 40/82, 40/84, 40/85, 40/89, 40/91, 40/92, 40/96, 40/98, 41/20, 41/21, 41/22, 41/23, 41/24, 41/27, 41/28, 41/29, 41/30, 41/31, 41/34, 41/35, 41/36, 41/41, 41/42, 41/43, 41/48, 41/49, 41/50, 41/55, 41/56, 41/57, 41/62, 41/63, 41/64, 41/65, 41/66, 41/67, 41/68, 41/71, 41/72, 41/73, 41/74, 41/75, 41/78, 41/79, 41/80, 41/81, 41/82, 41/85, 41/86, 41/87, 41/88, 41/89, 41/92, 41/93, 41/94, 41/95, 41/96, 42/21, 42/23, 42/24, 42/28, 42/30, 42/31, 42/35, 42/36, 42/42, 42/43, 42/49, 42/50, 42/56, 42/57, 42/63, 42/65, 42/66, 42/68, 42/72, 42/74, 42/75, 42/79, 42/81, 42/82, 42/86, 42/88, 42/89, 42/93, 42/95, 42/96, 43/20, 43/21, 43/24, 43/25, 43/26, 43/27, 43/28, 43/31, 43/36, 43/43, 43/50, 43/57, 43/58, 43/59, 43/60, 43/61, 43/62, 43/63, 43/66, 43/67, 43/68, 43/69, 43/70, 43/71, 43/72, 43/75, 43/76, 43/77, 43/78, 43/79, 43/82, 43/83, 43/84, 43/85, 43/86, 43/89, 43/90, 43/91, 43/92, 43/93, 43/96, 43/97, 43/98, 44/20, 44/21, 44/25, 44/26, 44/28, 44/29, 44/30, 44/31, 44/32, 44/33, 44/34, 44/35, 44/36, 44/37, 44/38, 44/39, 44/40, 44/41, 44/42, 44/43, 44/44, 44/45, 44/46, 44/47, 44/48, 44/49, 44/50, 44/58, 44/60, 44/62, 44/63, 44/67, 44/69, 44/70, 44/72, 44/76, 44/77, 44/79, 44/83, 44/84, 44/86, 44/90, 44/91, 44/93, 44/97, 44/98, 45/21, 45/26, 45/29, 45/30, 45/32, 45/33, 45/35, 45/37, 45/39, 45/40, 45/42, 45/44, 45/46, 45/48, 45/50, 45/51, 45/52, 45/53, 45/54, 45/55, 45/56, 45/57, 45/58, 45/59, 45/60, 45/63, 45/64, 45/65, 45/66, 45/67, 45/70, 45/71, 45/72, 45/77, 45/78, 45/79, 45/84, 45/85, 45/86, 45/91, 45/92, 45/93, 45/98, 46/20, 46/21, 46/22, 46/23, 46/24, 46/25, 46/26, 46/30, 46/33, 46/34, 46/35, 46/36, 46/37, 46/40, 46/41, 46/42, 46/43, 46/44, 46/45, 46/46, 46/47, 46/48, 46/51, 46/52, 46/54, 46/55, 46/57, 46/59, 46/60, 46/64, 46/66, 46/67, 46/71, 46/72, 46/78, 46/79, 46/85, 46/86, 46/92, 46/93, 47/20, 47/22, 47/24, 47/26, 47/27, 47/28, 47/29, 47/30, 47/34, 47/36, 47/37, 47/41, 47/43, 47/45, 47/47, 47/48, 47/52, 47/55, 47/56, 47/57, 47/60, 47/61, 47/62, 47/63, 47/64, 47/67, 47/72, 47/79, 47/86, 47/93, 47/94, 47/95, 47/96, 47/97, 47/98, 48/20, 48/21, 48/22, 48/23, 48/24, 48/27, 48/28, 48/30, 48/31, 48/32, 48/33, 48/34, 48/37, 48/38, 48/39, 48/40, 48/41, 48/42, 48/43, 48/44, 48/45, 48/48, 48/49, 48/50, 48/51, 48/52, 48/56, 48/57, 48/61, 48/62, 48/64, 48/65, 48/66, 48/67, 48/68, 48/69, 48/70, 48/71, 48/72, 48/73, 48/74, 48/75, 48/76, 48/77, 48/78, 48/79, 48/80, 48/81, 48/82, 48/83, 48/84, 48/85, 48/86, 48/94, 48/96, 48/98, 49/21, 49/23, 49/24, 49/28, 49/31, 49/32, 49/34, 49/38, 49/40, 49/42, 49/44, 49/45, 49/49, 49/51, 49/52, 49/57, 49/62, 49/65, 49/66, 49/68, 49/69, 49/71, 49/73, 49/75, 49/76, 49/78, 49/80, 49/82, 49/84, 49/86, 49/87, 49/88, 49/89, 49/90, 49/91, 49/92, 49/93, 49/94, 49/95, 49/96, 50/20, 50/21, 50/24, 50/25, 50/26, 50/27, 50/28, 50/32, 50/33, 50/34, 50/35, 50/36, 50/37, 50/38, 50/39, 50/40, 50/41, 50/42, 50/45, 50/46, 50/47, 50/48, 50/49, 50/52, 50/53, 50/54, 50/55, 50/56, 50/57, 50/58, 50/59, 50/60, 50/61, 50/62, 50/66, 50/69, 50/70, 50/71, 50/72, 50/73, 50/76, 50/77, 50/78, 50/79, 50/80, 50/81, 50/82, 50/83, 50/84, 50/87, 50/88, 50/90, 50/91, 50/93, 50/95, 50/96, 51/20, 51/21, 51/25, 51/27, 51/28, 51/33, 51/35, 51/37, 51/39, 51/41, 51/42, 51/46, 51/48, 51/49, 51/53, 51/55, 51/56, 51/58, 51/60, 51/62, 51/63, 51/64, 51/65, 51/66, 51/70, 51/72, 51/73, 51/77, 51/79, 51/81, 51/83, 51/84, 51/88, 51/91, 51/92, 51/93, 51/96, 51/97, 51/98, 52/21, 52/22, 52/23, 52/24, 52/25, 52/28, 52/29, 52/30, 52/31, 52/32, 52/33, 52/34, 52/35, 52/36, 52/37, 52/38, 52/39, 52/42, 52/43, 52/44, 52/45, 52/46, 52/49, 52/50, 52/51, 52/52, 52/53, 52/56, 52/57, 52/58, 52/59, 52/60, 52/63, 52/64, 52/66, 52/67, 52/68, 52/69, 52/70, 52/73, 52/74, 52/75, 52/76, 52/77, 52/78, 52/79, 52/80, 52/81, 52/84, 52/85, 52/86, 52/87, 52/88, 52/92, 52/93, 52/97, 52/98, 53/22, 53/24, 53/25, 53/29, 53/31, 53/32, 53/34, 53/36, 53/38, 53/39, 53/43, 53/45, 53/46, 53/50, 53/52, 53/53, 53/57, 53/59, 53/60, 53/64, 53/67, 53/68, 53/70, 53/74, 53/76, 53/78, 53/80, 53/81, 53/85, 53/87, 53/88, 53/93, 53/98, 54/20, 54/21, 54/22, 54/25, 54/26, 54/27, 54/28, 54/29, 54/32, 54/33, 54/34, 54/35, 54/36, 54/39, 54/40, 54/41, 54/42, 54/43, 54/46, 54/47, 54/48, 54/49, 54/50, 54/53, 54/54, 54/55, 54/56, 54/57, 54/60, 54/61, 54/62, 54/63, 54/64, 54/68, 54/69, 54/70, 54/71, 54/72, 54/73, 54/74, 54/75, 54/76, 54/77, 54/78, 54/81, 54/82, 54/83, 54/84, 54/85, 54/88, 54/89, 54/90, 54/91, 54/92, 54/93, 54/94, 54/95, 54/96, 54/97, 54/98, 55/21, 55/22, 55/26, 55/28, 55/29, 55/33, 55/35, 55/36, 55/40, 55/42, 55/43, 55/47, 55/49, 55/50, 55/54, 55/56, 55/57, 55/61, 55/63, 55/64, 55/69, 55/71, 55/73, 55/75, 55/77, 55/78, 55/82, 55/84, 55/85, 55/89, 55/91, 55/92, 55/94, 55/96, 55/98, 56/22, 56/23, 56/24, 56/25, 56/26, 56/29, 56/30, 56/31, 56/32, 56/33, 56/36, 56/37, 56/38, 56/39, 56/40, 56/43, 56/44, 56/45, 56/46, 56/47, 56/50, 56/51, 56/52, 56/53, 56/54, 56/57, 56/58, 56/59, 56/60, 56/61, 56/64, 56/65, 56/66, 56/67, 56/68, 56/69, 56/70, 56/71, 56/72, 56/73, 56/74, 56/75, 56/78, 56/79, 56/80, 56/81, 56/82, 56/85, 56/86, 56/87, 56/88, 56/89, 56/92, 56/93, 56/94, 56/95, 56/96, 57/23, 57/25, 57/26, 57/30, 57/32, 57/33, 57/37, 57/39, 57/40, 57/44, 57/46, 57/47, 57/51, 57/53, 57/54, 57/58, 57/60, 57/61, 57/65, 57/67, 57/68, 57/70, 57/72, 57/74, 57/75, 57/79, 57/81, 57/82, 57/86, 57/88, 57/89, 57/93, 57/95, 57/96, 58/20, 58/21, 58/22, 58/23, 58/26, 58/27, 58/28, 58/29, 58/30, 58/33, 58/34, 58/35, 58/36, 58/37, 58/40, 58/41, 58/42, 58/43, 58/44, 58/47, 58/48, 58/49, 58/50, 58/51, 58/54, 58/55, 58/56, 58/57, 58/58, 58/61, 58/62, 58/63, 58/64, 58/65, 58/68, 58/69, 58/70, 58/71, 58/72, 58/75, 58/76, 58/77, 58/78, 58/79, 58/82, 58/83, 58/84, 58/85, 58/86, 58/89, 58/90, 58/91, 58/92, 58/93, 58/96, 58/97, 58/98, 59/20, 59/22, 59/23, 59/27, 59/29, 59/30, 59/34, 59/36, 59/37, 59/41, 59/43, 59/44, 59/48, 59/50, 59/51, 59/55, 59/57, 59/58, 59/62, 59/64, 59/65, 59/69, 59/71, 59/72, 59/76, 59/78, 59/79, 59/83, 59/85, 59/86, 59/90, 59/92, 59/93, 59/97, 60/20, 60/23, 60/24, 60/25, 60/26, 60/27, 60/30, 60/31, 60/32, 60/33, 60/34, 60/37, 60/38, 60/39, 60/40, 60/41, 60/44, 60/45, 60/46, 60/47, 60/48, 60/51, 60/52, 60/53, 60/54, 60/55, 60/58, 60/59, 60/60, 60/61, 60/62, 60/65, 60/66, 60/67, 60/68, 60/69, 60/72, 60/73, 60/74, 60/75, 60/76, 60/79, 60/80, 60/81, 60/82, 60/83, 60/86, 60/87, 60/88, 60/89, 60/90, 60/93, 60/94, 60/95, 60/96, 60/97, 61/20, 61/24, 61/26, 61/27, 61/31, 61/33, 61/34, 61/38, 61/40, 61/41, 61/45, 61/47, 61/48, 61/52, 61/54, 61/55, 61/59, 61/61, 61/62, 61/66, 61/68, 61/69, 61/73, 61/75, 61/76, 61/80, 61/82, 61/83, 61/87, 61/89, 61/90, 61/94, 61/96, 61/97, 62/20, 62/21, 62/22, 62/23, 62/24, 62/27, 62/28, 62/29, 62/30, 62/31, 62/34, 62/35, 62/36, 62/37, 62/38, 62/41, 62/42, 62/43, 62/44, 62/45, 62/48, 62/49, 62/50, 62/51, 62/52, 62/55, 62/56, 62/57, 62/58, 62/59, 62/62, 62/63, 62/64, 62/65, 62/66, 62/69, 62/70, 62/71, 62/72, 62/73, 62/76, 62/77, 62/78, 62/79, 62/80, 62/83, 62/84, 62/85, 62/86, 62/87, 62/90, 62/91, 62/92, 62/93, 62/94, 62/97, 62/98, 63/21, 63/23, 63/24, 63/28, 63/30, 63/31, 63/35, 63/37, 63/38, 63/42, 63/44, 63/45, 63/49, 63/51, 63/52, 63/56, 63/58, 63/59, 63/63, 63/65, 63/66, 63/70, 63/72, 63/73, 63/77, 63/79, 63/80, 63/84, 63/86, 63/87, 63/91, 63/93, 63/94, 63/98, 64/20, 64/21, 64/24, 64/25, 64/26, 64/27, 64/28, 64/31, 64/32, 64/33, 64/34, 64/35, 64/38, 64/39, 64/40, 64/41, 64/42, 64/45, 64/46, 64/47, 64/48, 64/49, 64/52, 64/53, 64/54, 64/55, 64/56, 64/59, 64/60, 64/61, 64/62, 64/63, 64/66, 64/67, 64/68, 64/69, 64/70, 64/73, 64/74, 64/75, 64/76, 64/77, 64/80, 64/81, 64/82, 64/83, 64/84, 64/87, 64/88, 64/89, 64/90, 64/91, 64/94, 64/95, 64/96, 64/97, 64/98, 65/20, 65/21, 65/25, 65/27, 65/28, 65/32, 65/34, 65/35, 65/39, 65/41, 65/42, 65/46, 65/48, 65/49, 65/53, 65/55, 65/56, 65/60, 65/62, 65/63, 65/67, 65/69, 65/70, 65/74, 65/76, 65/77, 65/81, 65/83, 65/84, 65/88, 65/90, 65/91, 65/95, 65/97, 65/98, 66/21, 66/22, 66/23, 66/24, 66/25, 66/28, 66/29, 66/30, 66/31, 66/32, 66/35, 66/36, 66/37, 66/38, 66/39, 66/42, 66/43, 66/44, 66/45, 66/46, 66/49, 66/50, 66/51, 66/52, 66/53, 66/56, 66/57, 66/58, 66/59, 66/60, 66/63, 66/64, 66/65, 66/66, 66/67, 66/70, 66/71, 66/72, 66/73, 66/74, 66/77, 66/78, 66/79, 66/80, 66/81, 66/84, 66/85, 66/86, 66/87, 66/88, 66/91, 66/92, 66/93, 66/94, 66/95, 66/98, 67/22, 67/24, 67/25, 67/29, 67/31, 67/32, 67/36, 67/38, 67/39, 67/43, 67/45, 67/46, 67/50, 67/52, 67/53, 67/57, 67/59, 67/60, 67/64, 67/66, 67/67, 67/71, 67/73, 67/74, 67/78, 67/80, 67/81, 67/85, 67/87, 67/88, 67/92, 67/94, 67/95, 68/20, 68/21, 68/22, 68/25, 68/26, 68/27, 68/28, 68/29, 68/32, 68/33, 68/34, 68/35, 68/36, 68/39, 68/40, 68/41, 68/42, 68/43, 68/46, 68/47, 68/48, 68/49, 68/50, 68/53, 68/54, 68/55, 68/56, 68/57, 68/60, 68/61, 68/62, 68/63, 68/64, 68/67, 68/68, 68/69, 68/70, 68/71, 68/74, 68/75, 68/76, 68/77, 68/78, 68/81, 68/82, 68/83, 68/84, 68/85, 68/88, 68/89, 68/90, 68/91, 68/92, 68/95, 68/96, 68/97, 68/98, 69/21, 69/22, 69/26, 69/28, 69/29, 69/33, 69/35, 69/36, 69/40, 69/42, 69/43, 69/47, 69/49, 69/50, 69/54, 69/56, 69/57, 69/61, 69/63, 69/64, 69/68, 69/70, 69/71, 69/75, 69/77, 69/78, 69/82, 69/84, 69/85, 69/89, 69/91, 69/92, 69/96, 69/98, 70/22, 70/23, 70/24, 70/25, 70/26, 70/29, 70/30, 70/31, 70/32, 70/33, 70/36, 70/37, 70/38, 70/39, 70/40, 70/43, 70/44, 70/45, 70/46, 70/47, 70/50, 70/51, 70/52, 70/53, 70/54, 70/57, 70/58, 70/59, 70/60, 70/61, 70/64, 70/65, 70/66, 70/67, 70/68, 70/71, 70/72, 70/73, 70/74, 70/75, 70/78, 70/79, 70/80, 70/81, 70/82, 70/85, 70/86, 70/87, 70/88, 70/89, 70/92, 70/93, 70/94, 70/95, 70/96, 71/23, 71/25, 71/26, 71/30, 71/32, 71/33, 71/37, 71/39, 71/40, 71/44, 71/46, 71/47, 71/51, 71/53, 71/54, 71/58, 71/60, 71/61, 71/65, 71/67, 71/68, 71/72, 71/74, 71/75, 71/79, 71/81, 71/82, 71/86, 71/88, 71/89, 71/93, 71/95, 71/96, 72/20, 72/21, 72/22, 72/23, 72/26, 72/27, 72/28, 72/29, 72/30, 72/33, 72/34, 72/35, 72/36, 72/37, 72/40, 72/41, 72/42, 72/43, 72/44, 72/47, 72/48, 72/49, 72/50, 72/51, 72/54, 72/55, 72/56, 72/57, 72/58, 72/61, 72/62, 72/63, 72/64, 72/65, 72/68, 72/69, 72/70, 72/71, 72/72, 72/75, 72/76, 72/77, 72/78, 72/79, 72/82, 72/83, 72/84, 72/85, 72/86, 72/89, 72/90, 72/91, 72/92, 72/93, 72/96, 72/97, 72/98, 73/20, 73/22, 73/23, 73/27, 73/29, 73/30, 73/34, 73/36, 73/37, 73/41, 73/43, 73/44, 73/48, 73/50, 73/51, 73/55, 73/57, 73/58, 73/62, 73/64, 73/65, 73/69, 73/71, 73/72, 73/76, 73/78, 73/79, 73/83, 73/85, 73/86, 73/90, 73/92, 73/93, 73/97, 74/20, 74/23, 74/24, 74/25, 74/26, 74/27, 74/30, 74/31, 74/32, 74/33, 74/34, 74/37, 74/38, 74/39, 74/40, 74/41, 74/44, 74/45, 74/46, 74/47, 74/48, 74/51, 74/52, 74/53, 74/54, 74/55, 74/58, 74/59, 74/60, 74/61, 74/62, 74/65, 74/66, 74/67, 74/68, 74/69, 74/72, 74/73, 74/74, 74/75, 74/76, 74/79, 74/80, 74/81, 74/82, 74/83, 74/86, 74/87, 74/88, 74/89, 74/90, 74/93, 74/94, 74/95, 74/96, 74/97, 75/20, 75/24, 75/26, 75/27, 75/31, 75/33, 75/34, 75/38, 75/40, 75/41, 75/45, 75/47, 75/48, 75/52, 75/54, 75/55, 75/59, 75/61, 75/62, 75/66, 75/68, 75/69, 75/73, 75/75, 75/76, 75/80, 75/82, 75/83, 75/87, 75/89, 75/90, 75/94, 75/96, 75/97, 76/20, 76/21, 76/22, 76/23, 76/24, 76/27, 76/28, 76/29, 76/30, 76/31, 76/34, 76/35, 76/36, 76/37, 76/38, 76/41, 76/42, 76/43, 76/44, 76/45, 76/48, 76/49, 76/50, 76/51, 76/52, 76/55, 76/56, 76/57, 76/58, 76/59, 76/62, 76/63, 76/64, 76/65, 76/66, 76/69, 76/70, 76/71, 76/72, 76/73, 76/76, 76/77, 76/78, 76/79, 76/80, 76/83, 76/84, 76/85, 76/86, 76/87, 76/90, 76/91, 76/92, 76/93, 76/94, 76/97, 76/98, 77/21, 77/23, 77/24, 77/28, 77/30, 77/31, 77/35, 77/37, 77/38, 77/42, 77/44, 77/45, 77/49, 77/51, 77/52, 77/56, 77/58, 77/59, 77/63, 77/65, 77/66, 77/70, 77/72, 77/73, 77/77, 77/79, 77/80, 77/84, 77/86, 77/87, 77/91, 77/93, 77/94, 77/98, 78/20, 78/21, 78/24, 78/25, 78/26, 78/27, 78/28, 78/31, 78/32, 78/33, 78/34, 78/35, 78/38, 78/39, 78/40, 78/41, 78/42, 78/45, 78/46, 78/47, 78/48, 78/49, 78/52, 78/53, 78/54, 78/55, 78/56, 78/59, 78/60, 78/61, 78/62, 78/63, 78/66, 78/67, 78/68, 78/69, 78/70, 78/73, 78/74, 78/75, 78/76, 78/77, 78/80, 78/81, 78/82, 78/83, 78/84, 78/87, 78/88, 78/89, 78/90, 78/91, 78/94, 78/95, 78/96, 78/97, 78/98
}
 \filldraw[color = red, opacity=0.9] (\x, \y-10) rectangle (\x+1, \y-10+1);

\foreach \x/\y in {
 10/20, 10/21, 10/23, 10/26, 10/27, 10/28, 10/30, 10/33, 10/34, 10/35, 10/37, 10/40, 10/41, 10/42, 10/44, 10/47, 10/48, 10/49, 10/51, 10/54, 10/55, 10/56, 10/58, 10/61, 10/62, 10/63, 10/65, 10/68, 10/69, 10/70, 10/72, 10/75, 10/76, 10/77, 10/79, 10/82, 10/83, 10/84, 10/86, 10/89, 10/90, 10/91, 10/93, 10/96, 10/97, 10/98, 11/23, 11/24, 11/30, 11/31, 11/37, 11/38, 11/44, 11/45, 11/51, 11/52, 11/58, 11/59, 11/65, 11/66, 11/72, 11/73, 11/79, 11/80, 11/86, 11/87, 11/93, 11/94, 12/20, 12/23, 12/24, 12/25, 12/27, 12/30, 12/31, 12/32, 12/34, 12/37, 12/38, 12/39, 12/41, 12/44, 12/45, 12/46, 12/48, 12/51, 12/52, 12/53, 12/55, 12/58, 12/59, 12/60, 12/62, 12/65, 12/66, 12/67, 12/69, 12/72, 12/73, 12/74, 12/76, 12/79, 12/80, 12/81, 12/83, 12/86, 12/87, 12/88, 12/90, 12/93, 12/94, 12/95, 12/97, 13/20, 13/21, 13/27, 13/28, 13/34, 13/35, 13/41, 13/42, 13/48, 13/49, 13/55, 13/56, 13/62, 13/63, 13/69, 13/70, 13/76, 13/77, 13/83, 13/84, 13/90, 13/91, 13/97, 13/98, 14/20, 14/21, 14/22, 14/24, 14/27, 14/28, 14/29, 14/31, 14/34, 14/35, 14/36, 14/38, 14/41, 14/42, 14/43, 14/45, 14/48, 14/49, 14/50, 14/52, 14/55, 14/56, 14/57, 14/59, 14/62, 14/63, 14/64, 14/66, 14/69, 14/70, 14/71, 14/73, 14/76, 14/77, 14/78, 14/80, 14/83, 14/84, 14/85, 14/87, 14/90, 14/91, 14/92, 14/94, 14/97, 14/98, 15/24, 15/25, 15/31, 15/32, 15/38, 15/39, 15/45, 15/46, 15/52, 15/53, 15/59, 15/60, 15/66, 15/67, 15/73, 15/74, 15/80, 15/81, 15/87, 15/88, 15/94, 15/95, 16/21, 16/24, 16/25, 16/26, 16/28, 16/31, 16/32, 16/33, 16/35, 16/38, 16/39, 16/40, 16/42, 16/45, 16/46, 16/47, 16/49, 16/52, 16/53, 16/54, 16/56, 16/59, 16/60, 16/61, 16/63, 16/66, 16/67, 16/68, 16/70, 16/73, 16/74, 16/75, 16/77, 16/80, 16/81, 16/82, 16/84, 16/87, 16/88, 16/89, 16/91, 16/94, 16/95, 16/96, 16/98, 17/21, 17/22, 17/28, 17/29, 17/35, 17/36, 17/42, 17/43, 17/49, 17/50, 17/56, 17/57, 17/63, 17/64, 17/70, 17/71, 17/77, 17/78, 17/84, 17/85, 17/91, 17/92, 17/98, 18/21, 18/22, 18/23, 18/25, 18/28, 18/29, 18/30, 18/32, 18/35, 18/36, 18/37, 18/39, 18/42, 18/43, 18/44, 18/46, 18/49, 18/50, 18/51, 18/53, 18/56, 18/57, 18/58, 18/60, 18/63, 18/64, 18/65, 18/67, 18/70, 18/71, 18/72, 18/74, 18/77, 18/78, 18/79, 18/81, 18/84, 18/85, 18/86, 18/88, 18/91, 18/92, 18/93, 18/95, 18/98, 19/25, 19/26, 19/32, 19/33, 19/39, 19/40, 19/46, 19/47, 19/53, 19/54, 19/60, 19/61, 19/67, 19/68, 19/74, 19/75, 19/81, 19/82, 19/88, 19/89, 19/95, 19/96, 20/20, 20/22, 20/25, 20/26, 20/27, 20/29, 20/32, 20/33, 20/34, 20/36, 20/39, 20/40, 20/41, 20/43, 20/46, 20/47, 20/48, 20/50, 20/53, 20/54, 20/55, 20/57, 20/60, 20/61, 20/62, 20/64, 20/67, 20/68, 20/69, 20/71, 20/74, 20/75, 20/76, 20/78, 20/81, 20/82, 20/83, 20/85, 20/88, 20/89, 20/90, 20/92, 20/95, 20/96, 20/97, 21/22, 21/23, 21/29, 21/30, 21/36, 21/37, 21/43, 21/44, 21/50, 21/51, 21/57, 21/58, 21/64, 21/65, 21/71, 21/72, 21/78, 21/79, 21/85, 21/86, 21/92, 21/93, 22/22, 22/23, 22/24, 22/26, 22/29, 22/30, 22/31, 22/33, 22/36, 22/37, 22/38, 22/40, 22/43, 22/44, 22/45, 22/47, 22/50, 22/51, 22/52, 22/54, 22/57, 22/58, 22/59, 22/61, 22/64, 22/65, 22/66, 22/68, 22/71, 22/72, 22/73, 22/75, 22/78, 22/79, 22/80, 22/82, 22/85, 22/86, 22/87, 22/89, 22/92, 22/93, 22/94, 22/96, 23/20, 23/26, 23/27, 23/33, 23/34, 23/40, 23/41, 23/47, 23/48, 23/54, 23/55, 23/61, 23/62, 23/68, 23/69, 23/75, 23/76, 23/82, 23/83, 23/89, 23/90, 23/96, 23/97, 24/20, 24/21, 24/23, 24/26, 24/27, 24/28, 24/30, 24/33, 24/34, 24/35, 24/37, 24/40, 24/41, 24/42, 24/44, 24/47, 24/48, 24/49, 24/51, 24/54, 24/55, 24/56, 24/58, 24/61, 24/62, 24/63, 24/65, 24/68, 24/69, 24/70, 24/72, 24/75, 24/76, 24/77, 24/79, 24/82, 24/83, 24/84, 24/86, 24/89, 24/90, 24/91, 24/93, 24/96, 24/97, 24/98, 25/23, 25/24, 25/30, 25/31, 25/37, 25/38, 25/44, 25/45, 25/51, 25/52, 25/58, 25/59, 25/65, 25/66, 25/72, 25/73, 25/79, 25/80, 25/86, 25/87, 25/93, 25/94, 26/20, 26/23, 26/24, 26/25, 26/27, 26/30, 26/31, 26/32, 26/34, 26/37, 26/38, 26/39, 26/41, 26/44, 26/45, 26/46, 26/48, 26/51, 26/52, 26/53, 26/55, 26/58, 26/59, 26/60, 26/62, 26/65, 26/66, 26/67, 26/69, 26/72, 26/73, 26/74, 26/76, 26/79, 26/80, 26/81, 26/83, 26/86, 26/87, 26/88, 26/90, 26/93, 26/94, 26/95, 26/97, 27/20, 27/21, 27/27, 27/28, 27/34, 27/35, 27/41, 27/42, 27/48, 27/49, 27/55, 27/56, 27/62, 27/63, 27/69, 27/70, 27/76, 27/77, 27/83, 27/84, 27/90, 27/91, 27/97, 27/98, 28/20, 28/21, 28/22, 28/24, 28/27, 28/28, 28/29, 28/31, 28/34, 28/35, 28/36, 28/38, 28/41, 28/42, 28/43, 28/45, 28/48, 28/49, 28/50, 28/52, 28/55, 28/56, 28/57, 28/59, 28/62, 28/63, 28/64, 28/66, 28/69, 28/70, 28/71, 28/73, 28/76, 28/77, 28/78, 28/80, 28/83, 28/84, 28/85, 28/87, 28/90, 28/91, 28/92, 28/94, 28/97, 28/98, 29/24, 29/25, 29/31, 29/32, 29/38, 29/39, 29/45, 29/46, 29/52, 29/53, 29/59, 29/60, 29/66, 29/67, 29/73, 29/74, 29/80, 29/81, 29/87, 29/88, 29/94, 29/95, 30/21, 30/24, 30/25, 30/26, 30/28, 30/31, 30/32, 30/33, 30/35, 30/38, 30/39, 30/40, 30/42, 30/45, 30/46, 30/47, 30/49, 30/52, 30/53, 30/54, 30/56, 30/59, 30/60, 30/61, 30/63, 30/66, 30/67, 30/68, 30/70, 30/73, 30/74, 30/75, 30/77, 30/80, 30/81, 30/82, 30/84, 30/87, 30/88, 30/89, 30/91, 30/94, 30/95, 30/96, 30/98, 31/21, 31/22, 31/28, 31/29, 31/35, 31/36, 31/42, 31/43, 31/49, 31/50, 31/56, 31/57, 31/63, 31/64, 31/70, 31/71, 31/77, 31/78, 31/84, 31/85, 31/91, 31/92, 31/98, 32/21, 32/22, 32/23, 32/25, 32/28, 32/29, 32/30, 32/32, 32/35, 32/36, 32/37, 32/39, 32/42, 32/43, 32/44, 32/46, 32/49, 32/50, 32/51, 32/53, 32/56, 32/57, 32/58, 32/60, 32/63, 32/64, 32/65, 32/67, 32/70, 32/71, 32/72, 32/74, 32/77, 32/78, 32/79, 32/81, 32/84, 32/85, 32/86, 32/88, 32/91, 32/92, 32/93, 32/95, 32/98, 33/25, 33/26, 33/32, 33/33, 33/39, 33/40, 33/46, 33/47, 33/53, 33/54, 33/60, 33/61, 33/67, 33/68, 33/74, 33/75, 33/81, 33/82, 33/88, 33/89, 33/95, 33/96, 34/20, 34/25, 34/26, 34/27, 34/29, 34/32, 34/33, 34/34, 34/36, 34/39, 34/40, 34/41, 34/43, 34/46, 34/47, 34/48, 34/50, 34/53, 34/54, 34/55, 34/57, 34/60, 34/61, 34/62, 34/64, 34/67, 34/68, 34/69, 34/71, 34/74, 34/75, 34/76, 34/78, 34/81, 34/82, 34/83, 34/85, 34/88, 34/89, 34/90, 34/92, 34/95, 34/96, 34/97, 35/29, 35/30, 35/36, 35/37, 35/43, 35/44, 35/50, 35/51, 35/57, 35/58, 35/64, 35/65, 35/71, 35/72, 35/78, 35/79, 35/85, 35/86, 35/92, 35/93, 36/29, 36/30, 36/31, 36/33, 36/36, 36/37, 36/38, 36/40, 36/43, 36/44, 36/45, 36/47, 36/50, 36/51, 36/52, 36/54, 36/57, 36/58, 36/59, 36/61, 36/64, 36/65, 36/66, 36/68, 36/71, 36/72, 36/73, 36/75, 36/78, 36/79, 36/80, 36/82, 36/85, 36/86, 36/87, 36/89, 36/92, 36/93, 36/94, 36/96, 37/33, 37/34, 37/40, 37/41, 37/47, 37/48, 37/54, 37/55, 37/61, 37/62, 37/68, 37/69, 37/75, 37/76, 37/82, 37/83, 37/89, 37/90, 37/96, 37/97, 38/33, 38/34, 38/35, 38/40, 38/41, 38/42, 38/47, 38/48, 38/49, 38/54, 38/55, 38/56, 38/61, 38/62, 38/63, 38/65, 38/68, 38/69, 38/70, 38/72, 38/75, 38/76, 38/77, 38/79, 38/82, 38/83, 38/84, 38/86, 38/89, 38/90, 38/91, 38/93, 38/96, 38/97, 38/98, 39/65, 39/66, 39/72, 39/73, 39/79, 39/80, 39/86, 39/87, 39/93, 39/94, 40/65, 40/66, 40/67, 40/69, 40/72, 40/73, 40/74, 40/76, 40/79, 40/80, 40/81, 40/83, 40/86, 40/87, 40/88, 40/90, 40/93, 40/94, 40/95, 40/97, 41/69, 41/70, 41/76, 41/77, 41/83, 41/84, 41/90, 41/91, 41/97, 41/98, 42/69, 42/70, 42/71, 42/76, 42/77, 42/78, 42/83, 42/84, 42/85, 42/90, 42/91, 42/92, 42/97, 42/98, 47/46, 48/46, 48/47, 49/20, 49/22, 49/43, 49/46, 49/47, 49/48, 49/50, 50/22, 50/23, 50/43, 50/44, 50/50, 50/51, 51/22, 51/23, 51/24, 51/26, 51/40, 51/43, 51/44, 51/45, 51/47, 51/50, 51/51, 51/52, 51/54, 51/82, 52/20, 52/26, 52/27, 52/40, 52/41, 52/47, 52/48, 52/54, 52/55, 52/82, 52/83, 53/20, 53/21, 53/23, 53/26, 53/27, 53/28, 53/30, 53/37, 53/40, 53/41, 53/42, 53/44, 53/47, 53/48, 53/49, 53/51, 53/54, 53/55, 53/56, 53/58, 53/79, 53/82, 53/83, 53/84, 53/86, 54/23, 54/24, 54/30, 54/31, 54/37, 54/38, 54/44, 54/45, 54/51, 54/52, 54/58, 54/59, 54/79, 54/80, 54/86, 54/87, 55/20, 55/23, 55/24, 55/25, 55/27, 55/30, 55/31, 55/32, 55/34, 55/37, 55/38, 55/39, 55/41, 55/44, 55/45, 55/46, 55/48, 55/51, 55/52, 55/53, 55/55, 55/58, 55/59, 55/60, 55/62, 55/76, 55/79, 55/80, 55/81, 55/83, 55/86, 55/87, 55/88, 55/90, 56/20, 56/21, 56/27, 56/28, 56/34, 56/35, 56/41, 56/42, 56/48, 56/49, 56/55, 56/56, 56/62, 56/63, 56/76, 56/77, 56/83, 56/84, 56/90, 56/91, 57/20, 57/21, 57/22, 57/24, 57/27, 57/28, 57/29, 57/31, 57/34, 57/35, 57/36, 57/38, 57/41, 57/42, 57/43, 57/45, 57/48, 57/49, 57/50, 57/52, 57/55, 57/56, 57/57, 57/59, 57/62, 57/63, 57/64, 57/66, 57/73, 57/76, 57/77, 57/78, 57/80, 57/83, 57/84, 57/85, 57/87, 57/90, 57/91, 57/92, 57/94, 58/24, 58/25, 58/31, 58/32, 58/38, 58/39, 58/45, 58/46, 58/52, 58/53, 58/59, 58/60, 58/66, 58/67, 58/73, 58/74, 58/80, 58/81, 58/87, 58/88, 58/94, 58/95, 59/21, 59/24, 59/25, 59/26, 59/28, 59/31, 59/32, 59/33, 59/35, 59/38, 59/39, 59/40, 59/42, 59/45, 59/46, 59/47, 59/49, 59/52, 59/53, 59/54, 59/56, 59/59, 59/60, 59/61, 59/63, 59/66, 59/67, 59/68, 59/70, 59/73, 59/74, 59/75, 59/77, 59/80, 59/81, 59/82, 59/84, 59/87, 59/88, 59/89, 59/91, 59/94, 59/95, 59/96, 59/98, 60/21, 60/22, 60/28, 60/29, 60/35, 60/36, 60/42, 60/43, 60/49, 60/50, 60/56, 60/57, 60/63, 60/64, 60/70, 60/71, 60/77, 60/78, 60/84, 60/85, 60/91, 60/92, 60/98, 61/21, 61/22, 61/23, 61/25, 61/28, 61/29, 61/30, 61/32, 61/35, 61/36, 61/37, 61/39, 61/42, 61/43, 61/44, 61/46, 61/49, 61/50, 61/51, 61/53, 61/56, 61/57, 61/58, 61/60, 61/63, 61/64, 61/65, 61/67, 61/70, 61/71, 61/72, 61/74, 61/77, 61/78, 61/79, 61/81, 61/84, 61/85, 61/86, 61/88, 61/91, 61/92, 61/93, 61/95, 61/98, 62/25, 62/26, 62/32, 62/33, 62/39, 62/40, 62/46, 62/47, 62/53, 62/54, 62/60, 62/61, 62/67, 62/68, 62/74, 62/75, 62/81, 62/82, 62/88, 62/89, 62/95, 62/96, 63/20, 63/22, 63/25, 63/26, 63/27, 63/29, 63/32, 63/33, 63/34, 63/36, 63/39, 63/40, 63/41, 63/43, 63/46, 63/47, 63/48, 63/50, 63/53, 63/54, 63/55, 63/57, 63/60, 63/61, 63/62, 63/64, 63/67, 63/68, 63/69, 63/71, 63/74, 63/75, 63/76, 63/78, 63/81, 63/82, 63/83, 63/85, 63/88, 63/89, 63/90, 63/92, 63/95, 63/96, 63/97, 64/22, 64/23, 64/29, 64/30, 64/36, 64/37, 64/43, 64/44, 64/50, 64/51, 64/57, 64/58, 64/64, 64/65, 64/71, 64/72, 64/78, 64/79, 64/85, 64/86, 64/92, 64/93, 65/22, 65/23, 65/24, 65/26, 65/29, 65/30, 65/31, 65/33, 65/36, 65/37, 65/38, 65/40, 65/43, 65/44, 65/45, 65/47, 65/50, 65/51, 65/52, 65/54, 65/57, 65/58, 65/59, 65/61, 65/64, 65/65, 65/66, 65/68, 65/71, 65/72, 65/73, 65/75, 65/78, 65/79, 65/80, 65/82, 65/85, 65/86, 65/87, 65/89, 65/92, 65/93, 65/94, 65/96, 66/20, 66/26, 66/27, 66/33, 66/34, 66/40, 66/41, 66/47, 66/48, 66/54, 66/55, 66/61, 66/62, 66/68, 66/69, 66/75, 66/76, 66/82, 66/83, 66/89, 66/90, 66/96, 66/97, 67/20, 67/21, 67/23, 67/26, 67/27, 67/28, 67/30, 67/33, 67/34, 67/35, 67/37, 67/40, 67/41, 67/42, 67/44, 67/47, 67/48, 67/49, 67/51, 67/54, 67/55, 67/56, 67/58, 67/61, 67/62, 67/63, 67/65, 67/68, 67/69, 67/70, 67/72, 67/75, 67/76, 67/77, 67/79, 67/82, 67/83, 67/84, 67/86, 67/89, 67/90, 67/91, 67/93, 67/96, 67/97, 67/98, 68/23, 68/24, 68/30, 68/31, 68/37, 68/38, 68/44, 68/45, 68/51, 68/52, 68/58, 68/59, 68/65, 68/66, 68/72, 68/73, 68/79, 68/80, 68/86, 68/87, 68/93, 68/94, 69/20, 69/23, 69/24, 69/25, 69/27, 69/30, 69/31, 69/32, 69/34, 69/37, 69/38, 69/39, 69/41, 69/44, 69/45, 69/46, 69/48, 69/51, 69/52, 69/53, 69/55, 69/58, 69/59, 69/60, 69/62, 69/65, 69/66, 69/67, 69/69, 69/72, 69/73, 69/74, 69/76, 69/79, 69/80, 69/81, 69/83, 69/86, 69/87, 69/88, 69/90, 69/93, 69/94, 69/95, 69/97, 70/20, 70/21, 70/27, 70/28, 70/34, 70/35, 70/41, 70/42, 70/48, 70/49, 70/55, 70/56, 70/62, 70/63, 70/69, 70/70, 70/76, 70/77, 70/83, 70/84, 70/90, 70/91, 70/97, 70/98, 71/20, 71/21, 71/22, 71/24, 71/27, 71/28, 71/29, 71/31, 71/34, 71/35, 71/36, 71/38, 71/41, 71/42, 71/43, 71/45, 71/48, 71/49, 71/50, 71/52, 71/55, 71/56, 71/57, 71/59, 71/62, 71/63, 71/64, 71/66, 71/69, 71/70, 71/71, 71/73, 71/76, 71/77, 71/78, 71/80, 71/83, 71/84, 71/85, 71/87, 71/90, 71/91, 71/92, 71/94, 71/97, 71/98, 72/24, 72/25, 72/31, 72/32, 72/38, 72/39, 72/45, 72/46, 72/52, 72/53, 72/59, 72/60, 72/66, 72/67, 72/73, 72/74, 72/80, 72/81, 72/87, 72/88, 72/94, 72/95, 73/21, 73/24, 73/25, 73/26, 73/28, 73/31, 73/32, 73/33, 73/35, 73/38, 73/39, 73/40, 73/42, 73/45, 73/46, 73/47, 73/49, 73/52, 73/53, 73/54, 73/56, 73/59, 73/60, 73/61, 73/63, 73/66, 73/67, 73/68, 73/70, 73/73, 73/74, 73/75, 73/77, 73/80, 73/81, 73/82, 73/84, 73/87, 73/88, 73/89, 73/91, 73/94, 73/95, 73/96, 73/98, 74/21, 74/22, 74/28, 74/29, 74/35, 74/36, 74/42, 74/43, 74/49, 74/50, 74/56, 74/57, 74/63, 74/64, 74/70, 74/71, 74/77, 74/78, 74/84, 74/85, 74/91, 74/92, 74/98, 75/21, 75/22, 75/23, 75/25, 75/28, 75/29, 75/30, 75/32, 75/35, 75/36, 75/37, 75/39, 75/42, 75/43, 75/44, 75/46, 75/49, 75/50, 75/51, 75/53, 75/56, 75/57, 75/58, 75/60, 75/63, 75/64, 75/65, 75/67, 75/70, 75/71, 75/72, 75/74, 75/77, 75/78, 75/79, 75/81, 75/84, 75/85, 75/86, 75/88, 75/91, 75/92, 75/93, 75/95, 75/98, 76/25, 76/26, 76/32, 76/33, 76/39, 76/40, 76/46, 76/47, 76/53, 76/54, 76/60, 76/61, 76/67, 76/68, 76/74, 76/75, 76/81, 76/82, 76/88, 76/89, 76/95, 76/96, 77/20, 77/22, 77/25, 77/26, 77/27, 77/29, 77/32, 77/33, 77/34, 77/36, 77/39, 77/40, 77/41, 77/43, 77/46, 77/47, 77/48, 77/50, 77/53, 77/54, 77/55, 77/57, 77/60, 77/61, 77/62, 77/64, 77/67, 77/68, 77/69, 77/71, 77/74, 77/75, 77/76, 77/78, 77/81, 77/82, 77/83, 77/85, 77/88, 77/89, 77/90, 77/92, 77/95, 77/96, 77/97, 78/22, 78/23, 78/29, 78/30, 78/36, 78/37, 78/43, 78/44, 78/50, 78/51, 78/57, 78/58, 78/64, 78/65, 78/71, 78/72, 78/78, 78/79, 78/85, 78/86, 78/92, 78/93
}
\filldraw[color = orange, opacity=0.4] (\x, \y-10) rectangle (\x+1, \y-10+1);

\draw[step=1cm,gray,ultra thin] (10, 10) grid (79,89); 

\draw (42, 0) node {Tape};
\draw (2, 55) node {T};\draw (2, 50) node {i};\draw (2,45) node {m};\draw (2, 40) node {e};

\end{tikzpicture}

%% file: CA2.tex
\begin{tikzpicture}[scale = 0.23]
\colorlet{darkgreen}{green!40!black}
\colorlet{lightgreen}{green!65!black}
\foreach \x/\y in {
 31/44, 31/45, 31/46, 31/47, 31/48, 31/51, 31/52, 31/53, 31/54, 31/55, 31/58, 31/59, 31/60, 31/61, 31/62, 31/65, 31/66, 31/67, 31/68, 31/69, 31/72, 31/73, 31/74, 31/75, 31/76, 31/79, 31/80, 31/81, 31/82, 31/83, 31/86, 31/87, 32/45, 32/47, 32/48, 32/52, 32/54, 32/55, 32/59, 32/61, 32/62, 32/66, 32/68, 32/69, 32/73, 32/75, 32/76, 32/80, 32/82, 32/83, 32/87, 33/43, 33/44, 33/45, 33/48, 33/49, 33/50, 33/51, 33/52, 33/55, 33/56, 33/57, 33/58, 33/59, 33/62, 33/63, 33/64, 33/65, 33/66, 33/69, 33/70, 33/71, 33/72, 33/73, 33/76, 33/77, 33/78, 33/79, 33/80, 33/83, 33/84, 33/85, 33/86, 33/87, 34/44, 34/45, 34/49, 34/51, 34/52, 34/56, 34/58, 34/59, 34/63, 34/65, 34/66, 34/70, 34/72, 34/73, 34/77, 34/79, 34/80, 34/84, 34/86, 34/87, 35/45, 35/46, 35/47, 35/48, 35/49, 35/52, 35/53, 35/54, 35/55, 35/56, 35/59, 35/60, 35/61, 35/62, 35/63, 35/66, 35/67, 35/68, 35/69, 35/70, 35/73, 35/74, 35/75, 35/76, 35/77, 35/80, 35/81, 35/82, 35/83, 35/84, 35/87, 36/46, 36/48, 36/49, 36/53, 36/55, 36/56, 36/60, 36/62, 36/63, 36/67, 36/69, 36/70, 36/74, 36/76, 36/77, 36/81, 36/83, 36/84, 37/43, 37/44, 37/45, 37/46, 37/49, 37/50, 37/51, 37/52, 37/53, 37/56, 37/57, 37/58, 37/59, 37/60, 37/63, 37/64, 37/65, 37/66, 37/67, 37/70, 37/71, 37/72, 37/73, 37/74, 37/77, 37/78, 37/79, 37/80, 37/81, 37/84, 37/85, 37/86, 37/87, 38/43, 38/45, 38/46, 38/50, 38/52, 38/53, 38/57, 38/59, 38/60, 38/64, 38/66, 38/67, 38/71, 38/73, 38/74, 38/78, 38/80, 38/81, 38/85, 38/87, 39/43, 39/46, 39/47, 39/48, 39/49, 39/50, 39/53, 39/54, 39/55, 39/56, 39/57, 39/60, 39/61, 39/62, 39/63, 39/64, 39/67, 39/68, 39/69, 39/70, 39/71, 39/74, 39/75, 39/76, 39/77, 39/78, 39/81, 39/82, 39/83, 39/84, 39/85, 40/43, 40/47, 40/48, 40/50, 40/54, 40/55, 40/57, 40/61, 40/62, 40/64, 40/68, 40/70, 40/71, 40/75, 40/77, 40/78, 40/82, 40/84, 40/85, 41/43, 41/48, 41/49, 41/50, 41/55, 41/56, 41/57, 41/62, 41/63, 41/64, 41/65, 41/66, 41/67, 41/68, 41/71, 41/72, 41/73, 41/74, 41/75, 41/78, 41/79, 41/80, 41/81, 41/82, 41/85, 41/86, 41/87, 42/43, 42/49, 42/50, 42/56, 42/57, 42/63, 42/65, 42/66, 42/68, 42/72, 42/74, 42/75, 42/79, 42/81, 42/82, 42/86, 43/43, 43/50, 43/57, 43/58, 43/59, 43/60, 43/61, 43/62, 43/63, 43/66, 43/67, 43/68, 43/69, 43/70, 43/71, 43/72, 43/75, 43/76, 43/77, 43/78, 43/79, 43/82, 43/83, 43/84, 43/85, 43/86, 44/43, 44/44, 44/45, 44/46, 44/47, 44/48, 44/49, 44/50, 44/58, 44/60, 44/62, 44/63, 44/67, 44/69, 44/70, 44/72, 44/76, 44/77, 44/79, 44/83, 44/84, 44/86, 45/44, 45/46, 45/48, 45/50, 45/51, 45/52, 45/53, 45/54, 45/55, 45/56, 45/57, 45/58, 45/59, 45/60, 45/63, 45/64, 45/65, 45/66, 45/67, 45/70, 45/71, 45/72, 45/77, 45/78, 45/79, 45/84, 45/85, 45/86, 46/43, 46/44, 46/45, 46/46, 46/47, 46/48, 46/51, 46/52, 46/54, 46/55, 46/57, 46/59, 46/60, 46/64, 46/66, 46/67, 46/71, 46/72, 46/78, 46/79, 46/85, 46/86, 47/43, 47/45, 47/47, 47/48, 47/52, 47/55, 47/56, 47/57, 47/60, 47/61, 47/62, 47/63, 47/64, 47/67, 47/72, 47/79, 47/86, 48/43, 48/44, 48/45, 48/48, 48/49, 48/50, 48/51, 48/52, 48/56, 48/57, 48/61, 48/62, 48/64, 48/65, 48/66, 48/67, 48/68, 48/69, 48/70, 48/71, 48/72, 48/73, 48/74, 48/75, 48/76, 48/77, 48/78, 48/79, 48/80, 48/81, 48/82, 48/83, 48/84, 48/85, 48/86, 49/44, 49/45, 49/49, 49/51, 49/52, 49/57, 49/62, 49/65, 49/66, 49/68, 49/69, 49/71, 49/73, 49/75, 49/76, 49/78, 49/80, 49/82, 49/84, 49/86, 49/87, 50/45, 50/46, 50/47, 50/48, 50/49, 50/52, 50/53, 50/54, 50/55, 50/56, 50/57, 50/58, 50/59, 50/60, 50/61, 50/62, 50/66, 50/69, 50/70, 50/71, 50/72, 50/73, 50/76, 50/77, 50/78, 50/79, 50/80, 50/81, 50/82, 50/83, 50/84, 50/87, 51/46, 51/48, 51/49, 51/53, 51/55, 51/56, 51/58, 51/60, 51/62, 51/63, 51/64, 51/65, 51/66, 51/70, 51/72, 51/73, 51/77, 51/79, 51/81, 51/83, 51/84, 52/43, 52/44, 52/45, 52/46, 52/49, 52/50, 52/51, 52/52, 52/53, 52/56, 52/57, 52/58, 52/59, 52/60, 52/63, 52/64, 52/66, 52/67, 52/68, 52/69, 52/70, 52/73, 52/74, 52/75, 52/76, 52/77, 52/78, 52/79, 52/80, 52/81, 52/84, 52/85, 52/86, 52/87, 53/43, 53/45, 53/46, 53/50, 53/52, 53/53, 53/57, 53/59, 53/60, 53/64, 53/67, 53/68, 53/70, 53/74, 53/76, 53/78, 53/80, 53/81, 53/85, 53/87, 54/43, 54/46, 54/47, 54/48, 54/49, 54/50, 54/53, 54/54, 54/55, 54/56, 54/57, 54/60, 54/61, 54/62, 54/63, 54/64, 54/68, 54/69, 54/70, 54/71, 54/72, 54/73, 54/74, 54/75, 54/76, 54/77, 54/78, 54/81, 54/82, 54/83, 54/84, 54/85, 55/43, 55/47, 55/49, 55/50, 55/54, 55/56, 55/57, 55/61, 55/63, 55/64, 55/69, 55/71, 55/73, 55/75, 55/77, 55/78, 55/82, 55/84, 55/85, 56/43, 56/44, 56/45, 56/46, 56/47, 56/50, 56/51, 56/52, 56/53, 56/54, 56/57, 56/58, 56/59, 56/60, 56/61, 56/64, 56/65, 56/66, 56/67, 56/68, 56/69, 56/70, 56/71, 56/72, 56/73, 56/74, 56/75, 56/78, 56/79, 56/80, 56/81, 56/82, 56/85, 56/86, 56/87, 57/44, 57/46, 57/47, 57/51, 57/53, 57/54, 57/58, 57/60, 57/61, 57/65, 57/67, 57/68, 57/70, 57/72, 57/74, 57/75, 57/79, 57/81, 57/82, 57/86, 58/43, 58/44, 58/47, 58/48, 58/49, 58/50, 58/51, 58/54, 58/55, 58/56, 58/57, 58/58, 58/61, 58/62, 58/63, 58/64, 58/65, 58/68, 58/69, 58/70, 58/71, 58/72, 58/75, 58/76, 58/77, 58/78, 58/79, 58/82, 58/83, 58/84, 58/85, 58/86, 59/43, 59/44, 59/48, 59/50, 59/51, 59/55, 59/57, 59/58, 59/62, 59/64, 59/65, 59/69, 59/71, 59/72, 59/76, 59/78, 59/79, 59/83, 59/85, 59/86, 60/44, 60/45, 60/46, 60/47, 60/48, 60/51, 60/52, 60/53, 60/54, 60/55, 60/58, 60/59, 60/60, 60/61, 60/62, 60/65, 60/66, 60/67, 60/68, 60/69, 60/72, 60/73, 60/74, 60/75, 60/76, 60/79, 60/80, 60/81, 60/82, 60/83, 60/86, 60/87, 61/45, 61/47, 61/48, 61/52, 61/54, 61/55, 61/59, 61/61, 61/62, 61/66, 61/68, 61/69, 61/73, 61/75, 61/76, 61/80, 61/82, 61/83, 61/87, 62/43, 62/44, 62/45, 62/48, 62/49, 62/50, 62/51, 62/52, 62/55, 62/56, 62/57, 62/58, 62/59, 62/62, 62/63, 62/64, 62/65, 62/66, 62/69, 62/70, 62/71, 62/72, 62/73, 62/76, 62/77, 62/78, 62/79, 62/80, 62/83, 62/84, 62/85, 62/86, 62/87, 63/44, 63/45, 63/49, 63/51, 63/52, 63/56, 63/58, 63/59, 63/63, 63/65, 63/66, 63/70, 63/72, 63/73, 63/77, 63/79, 63/80, 63/84, 63/86, 63/87
}
\filldraw[color = red, opacity=0.7] (\x-18, \y-36) rectangle (\x-18+1, \y-36+1);

\foreach \x/\y in {
 31/43, 31/49, 31/50, 31/56, 31/57, 31/63, 31/64, 31/70, 31/71, 31/77, 31/78, 31/84, 31/85, 32/43, 32/44, 32/46, 32/49, 32/50, 32/51, 32/53, 32/56, 32/57, 32/58, 32/60, 32/63, 32/64, 32/65, 32/67, 32/70, 32/71, 32/72, 32/74, 32/77, 32/78, 32/79, 32/81, 32/84, 32/85, 32/86, 33/46, 33/47, 33/53, 33/54, 33/60, 33/61, 33/67, 33/68, 33/74, 33/75, 33/81, 33/82, 34/43, 34/46, 34/47, 34/48, 34/50, 34/53, 34/54, 34/55, 34/57, 34/60, 34/61, 34/62, 34/64, 34/67, 34/68, 34/69, 34/71, 34/74, 34/75, 34/76, 34/78, 34/81, 34/82, 34/83, 34/85, 35/43, 35/44, 35/50, 35/51, 35/57, 35/58, 35/64, 35/65, 35/71, 35/72, 35/78, 35/79, 35/85, 35/86, 36/43, 36/44, 36/45, 36/47, 36/50, 36/51, 36/52, 36/54, 36/57, 36/58, 36/59, 36/61, 36/64, 36/65, 36/66, 36/68, 36/71, 36/72, 36/73, 36/75, 36/78, 36/79, 36/80, 36/82, 36/85, 36/86, 36/87, 37/47, 37/48, 37/54, 37/55, 37/61, 37/62, 37/68, 37/69, 37/75, 37/76, 37/82, 37/83, 38/47, 38/48, 38/49, 38/54, 38/55, 38/56, 38/61, 38/62, 38/63, 38/65, 38/68, 38/69, 38/70, 38/72, 38/75, 38/76, 38/77, 38/79, 38/82, 38/83, 38/84, 38/86, 39/65, 39/66, 39/72, 39/73, 39/79, 39/80, 39/86, 39/87, 40/65, 40/66, 40/67, 40/69, 40/72, 40/73, 40/74, 40/76, 40/79, 40/80, 40/81, 40/83, 40/86, 40/87, 41/69, 41/70, 41/76, 41/77, 41/83, 41/84, 42/69, 42/70, 42/71, 42/76, 42/77, 42/78, 42/83, 42/84, 42/85, 47/46, 48/46, 48/47, 49/43, 49/46, 49/47, 49/48, 49/50, 50/43, 50/44, 50/50, 50/51, 51/43, 51/44, 51/45, 51/47, 51/50, 51/51, 51/52, 51/54, 51/82, 52/47, 52/48, 52/54, 52/55, 52/82, 52/83, 53/44, 53/47, 53/48, 53/49, 53/51, 53/54, 53/55, 53/56, 53/58, 53/79, 53/82, 53/83, 53/84, 53/86, 54/44, 54/45, 54/51, 54/52, 54/58, 54/59, 54/79, 54/80, 54/86, 54/87, 55/44, 55/45, 55/46, 55/48, 55/51, 55/52, 55/53, 55/55, 55/58, 55/59, 55/60, 55/62, 55/76, 55/79, 55/80, 55/81, 55/83, 55/86, 55/87, 56/48, 56/49, 56/55, 56/56, 56/62, 56/63, 56/76, 56/77, 56/83, 56/84, 57/43, 57/45, 57/48, 57/49, 57/50, 57/52, 57/55, 57/56, 57/57, 57/59, 57/62, 57/63, 57/64, 57/66, 57/73, 57/76, 57/77, 57/78, 57/80, 57/83, 57/84, 57/85, 57/87, 58/45, 58/46, 58/52, 58/53, 58/59, 58/60, 58/66, 58/67, 58/73, 58/74, 58/80, 58/81, 58/87, 59/45, 59/46, 59/47, 59/49, 59/52, 59/53, 59/54, 59/56, 59/59, 59/60, 59/61, 59/63, 59/66, 59/67, 59/68, 59/70, 59/73, 59/74, 59/75, 59/77, 59/80, 59/81, 59/82, 59/84, 59/87, 60/43, 60/49, 60/50, 60/56, 60/57, 60/63, 60/64, 60/70, 60/71, 60/77, 60/78, 60/84, 60/85, 61/43, 61/44, 61/46, 61/49, 61/50, 61/51, 61/53, 61/56, 61/57, 61/58, 61/60, 61/63, 61/64, 61/65, 61/67, 61/70, 61/71, 61/72, 61/74, 61/77, 61/78, 61/79, 61/81, 61/84, 61/85, 61/86, 62/46, 62/47, 62/53, 62/54, 62/60, 62/61, 62/67, 62/68, 62/74, 62/75, 62/81, 62/82, 63/43, 63/46, 63/47, 63/48, 63/50, 63/53, 63/54, 63/55, 63/57, 63/60, 63/61, 63/62, 63/64, 63/67, 63/68, 63/69, 63/71, 63/74, 63/75, 63/76, 63/78, 63/81, 63/82, 63/83, 63/85
}
\filldraw[color = orange, opacity=0.6] (\x-18, \y-36) rectangle (\x-18+1, \y-36+1);

\foreach \x/\y in {
20/14,
26/22,
26/24,
28/28,
30/30,
33/34,
35/38,
35/40,
35/42
}
 \filldraw[color = darkgreen, opacity=1]  (\x+.5, \y+.5) circle (10pt);
\foreach \x/\y in {
20/15, 21/15, 22/15, 23/15, 24/15, 25/15,
26/21, 26/20, 26/19, 26/18, 26/17, 26/16, 26/15,
26/23,
26/25, 27/25, 28/25, 28/26, 28/27, 
29/29,
31/31, 32/31, 33/31, 
33/32,
33/33,
33/35, 34/35, 35/35, 
35/36,
35/37,
35/39,
35/41,
35/43,
36/43,
37/43,
37/44,
37/45
}
\filldraw[color = lightgreen, opacity=1]  (\x+.5, \y+.5) circle (10pt);
\draw[step=1cm,gray,ultra thin] (13, 7) grid (46,52); 

\draw (29, 2) node {Tape};
\draw (7, 35) node {T};\draw (7, 33) node {i};\draw (7,31) node {m};\draw (7, 29) node {e};

\end{tikzpicture}